\documentclass[
11pt]{article}
\usepackage{amssymb}
\usepackage{amsmath}
\usepackage{color}
\usepackage{mathrsfs}
\numberwithin{equation}{section}
\newtheorem{theo}{Theorem}[section]
\newtheorem{defi}[theo]{Definition}
\newtheorem{rem}[theo]{\it Remark}
\newtheorem{prop}[theo]{Proposition}

\newtheorem{lemma}[theo]{Lemma}

\def\pp{{\bf P}}
\def\ppe{\pp_{\hspace{-3pt}\varepsilon}}
\def\ppo{\pp_{\hspace{-3pt}{\gamma}}}
\def\ee{{\bf E}}
\def\eeo{\ee_{\hspace{-0pt}{\gamma}}}
\def\dd{{\rm d}}

\def\Box{\rule{6pt}{6pt}}

\def\CJ{C_{\hspace{-2pt}J}}
\def\NF{\mathbb{F}}

\def\NR{\mathbb{R}}

\def\NZ{\mathbb{Z}}
\begin{document}
\title{A rigorous equation for the Cole-Hopf solution of the\\
conservative KPZ equation}
\author{Sigurd Assing\\
Department of Statistics, The University of Warwick,
Coventry CV4 7AL, UK\\
e-mail: s.assing@warwick.ac.uk}
\date{}
\maketitle
\begin{abstract}
A rigorous equation is stated and it is shown that the spatial
derivative of the Cole-Hopf solution of the KPZ equation 
is a solution of this equation. 
The method of proof used to show that a process solves this equation
is based on rather weak estimates so that this method has the advantage that
it could be used to verify solutions of other highly singular SPDEs, too.
\end{abstract}
%
\noindent
{\large KEY WORDS}\hspace{0.4cm}
stochastic partial differential equation,
interacting particle system, 
martingale problem,
generalized functions

\vspace{10pt}
\noindent
{\it Mathematics Subject Classification (2000)}:
Primary 60H15; Secondary 60K35
\section{Motivation and Summary}\label{moti}
The formal equation discussed in this paper is
\begin{equation}\label{star}
\frac{\partial}{\partial t}Y
\,=\,
\frac{\partial^2}{\partial u^2}Y
+{\gamma}\,\frac{\partial}{\partial u}(Y^2)
+\sqrt{2}\,\frac{\partial^3 B}{\partial t\partial u^2}
\end{equation}
where ${\gamma}$ is a real-valued parameter
and $B$ stands for a Brownian sheet thus 
$\frac{\partial^3 B}{\partial t\partial u^2}$
can be interpreted as the spatial derivative of a 
space-time white noise driving force.
The potential solutions $Y$ to this equation which were first
constructed in \cite{BG1997} take values in the space
$D([0,T];{\mathscr D}^\prime(\NR))$
of all cadlag functions mapping $[0,T]$ into the space of
Schwartz distributions ${\mathscr D}^\prime(\NR)$.
So the problem arises to give meaning to the non-linear term
$\frac{\partial}{\partial u}(Y^2)$ and this is meant by
{\it stating a rigorous equation} in this paper.

Equation $(\ref{star})$ is the equation the spatial derivative of a
solution
of the KPZ equation for growing interfaces would formally satisfy and
the main result in \cite{BG1997} is actually an approximation scheme
for the KPZ equation. The limiting field of this approximation scheme
equals the Cole-Hopf transform of another process and the community
started to call it the Cole-Hopf solution of the KPZ equation.
Taking the spatial derivative of the KPZ
equation turns it into a conservative system with an invariant state
and that's why $(\ref{star})$ is also called conservative KPZ
equation.

There has been a recent breakthrough in the theory of solutions to the
KPZ equation, see \cite{H2012}, and the reader is referred to this
work and the references therein for a good account on the progress
being made over the past few years in the understanding of the KPZ
equation. But, as in \cite{BG1997}, the main focus in \cite{H2012} is
on an approximation scheme and it is not shown that the limiting
field,
which again equals the Cole-Hopf solution, is the solution of a
well-defined equation.

Since $Y_t\in{\mathscr D}^\prime(\NR)$ for fixed $t$, 
the canonical definition
of the ill-posed term $\frac{\partial}{\partial u}(Y_t^2)$ 
would be a limit of type
$\frac{\partial}{\partial u}[(Y_t\star J_N)^2],\,N\to\infty$, 
using a mollifier $J\in{\mathscr D}(\NR)$ to approximate the
identity. Here $Y_t\star J_N$ denotes the convolution of the generalized function
$Y_t$ and the smooth function $J_N(u)=NJ(Nu),\,u\in\NR$.

It turned out that, 
even in the case where $Y_t$ is stationary,
it is hard to make sense of such a
limit in an appropriate space. The author only achieved to get
convergence in a rather artificial space of so-called generalized
random variables
which made it kind of impossible to understand $(\ref{star})$ as a PDE and
the notion of solution was based on a generalized martingale problem
(see \cite{A2002}). It even remains to be shown that $Y$
is indeed a solution of this generalized martingale problem.

The difficulty seems to be that, as far as we know,
there is no control of moments higher than two. 
Very good if not the best second order moment estimates for 
$Y_t(G)$ in the case where $Y_t$ is stationary
can be found in \cite{BQS2011} but the authors
themselves remark that their method cannot be applied to moments of
higher order. 
 
On the other hand, the convergence of time integrals 
$\int_r^t\frac{\partial}{\partial u}[(Y_s\star J_N)^2]\,\dd s,\,N\to\infty$, 
$r\le t$ fixed,
is much more regular and the notion of solution to (\ref{star})
introduced in \cite{JG2010} is based on the existence of such a limit. 
However, in \cite{JG2010} it is not explained how 
$\frac{\partial}{\partial u}(Y^2)$ should be understood for a chosen
$Y\in D([0,T];{\mathscr D}^\prime(\NR))$.
Instead, first showing very useful estimates, 
the authors of \cite{JG2010} conclude that
\begin{equation}\label{mean square}
-\lim_{N\to\infty}\int_r^t\int_\NR(Y_s\star J_N)^2(u)
\frac{G(u+1/N)-G(u)}{1/N}\,\dd u\dd s
\quad\mbox{exists in mean square}
\end{equation}
for every $r\le t$ 
and every test function $G$ in the Schwartz space ${\mathscr S}(\NR)$.
If $\frac{\partial}{\partial u}(Y^2)$ is defined by a limit
for every $Y\in D([0,T];{\mathscr D}^\prime(\NR))$ then verifying
equation $(\ref{star})$ for a possible solution requires a further
limit-exchange and this has not been accomplished in \cite{JG2010}.

The main message from \cite{A2002} is that
interchanging $\lim_{N\to\infty}$ and the time integration
in (\ref{mean square}) leads to severe complications.
So one wants to define
$$\langle{\bf 1}_{[r,t]}\otimes G\,,\,
\frac{\partial}{\partial u}(Y^2)\rangle
\quad\mbox{by}\quad
-\lim_{N\to\infty}
\int_r^t\int_\NR G^\prime(u)\,(Y_s\star J_N)^2(u)\,\dd u\dd s,$$
thinking of ${\bf 1}_{[r,t]}\otimes G$ as a test function and
of $\langle\,\,,\,\rangle$ as a dual pairing,
which triggers the idea to explain $\frac{\partial}{\partial u}(Y^2)$
as an element of ${\mathscr D}^\prime((0,T)\times\NR)$. 
Indeed, if
$$\lim_{N\to\infty}\int_0^T\hspace{-7pt}\int_\NR
\frac{\partial}{\partial u}\phi(t,u)\,(Y_t\star J_{N})^2(u)\,\dd u\dd t
\quad\mbox{exists for all
$\phi\in{\mathscr D}((0,T)\times\NR)$}$$
then this defines an element in ${\mathscr D}^\prime((0,T)\times\NR)$.
Of course, the above limit does not exist for all 
$Y\in D([0,T];{\mathscr D}^\prime(\NR))$ 
and limits of subsequences can be different depending on $\phi$.
So the definition of $\frac{\partial}{\partial u}(Y^2)$
justified in this paper requires finding a suitable subsequence
$(N_k)_{k=1}^\infty$ which is used to split
$D([0,T];{\mathscr D}^\prime(\NR))$ into two sets 
${\mathscr N}_{div}\cup{\mathscr N}_{div}^c$ where
\begin{eqnarray}
{\mathscr N}_{div}\,\stackrel{\mbox{\tiny def}}{=}\,
\left\{\rule{0pt}{13pt}\right.
Y\in D([0,T];{\mathscr D}^\prime(\NR))&:&
\lim_{k\to\infty}\int_0^T\hspace{-7pt}\int_\NR
\frac{\partial}{\partial u}\phi(t,u)
\,(Y_t\star J_{N_k})^2(u)\,\dd u\dd t \nonumber\\
&&\mbox{does not exist for some 
$\phi\in{\mathscr D}((0,T)\times\NR)$}
\left.\rule{0pt}{13pt}\right\}. \label{Ndiv}
\end{eqnarray}
Defining 
$\frac{\partial}{\partial u}(Y^2)\in{\mathscr D}^\prime((0,T)\times\NR)$
for every $\phi\in{\mathscr D}((0,T)\times\NR)$ by
\begin{equation}\label{first def}
\langle\phi\,,\,\frac{\partial}{\partial u}(Y^2)\rangle
\,\stackrel{\mbox{\tiny def}}{=}\left\{\begin{array}{rcl}
0&:&Y\in{\mathscr N}_{div}\\
\rule{0pt}{15pt}
-\lim_{k\to\infty}\int_0^T\hspace{-5pt}\int_\NR
\frac{\partial}{\partial u}\phi(t,u)
\,(Y_t\star J_{N_k})^2(u)\,\dd u\dd t&:&Y\in{\mathscr N}_{div}^c
\end{array}\right.
\end{equation}
turns the equation $(\ref{star})$ into a classical SPDE 
and it will be shown in this paper that
\begin{equation}\label{weak star}
\langle\phi\,,\,
\frac{\partial}{\partial t}Y
-\frac{\partial^2}{\partial u^2}Y
-{\gamma}\,\frac{\partial}{\partial u}(Y^2)
-\sqrt{2}\,\frac{\partial^3 B}{\partial t\partial u^2}\rangle
\,=\,0
\quad\mbox{for all $\phi\in{\mathscr D}((0,T)\times\NR)$ a.s.}
\end{equation}
for the stationary (potential) solution $Y$ constructed in \cite{BG1997} 
and some Brownian sheet $B$ both given on the same probability space.

Notice that the limits defining $\frac{\partial}{\partial u}(Y^2)$ in
the case where $Y\in{\mathscr N}_{div}^c$ could depend on the choice
of the mollifier $J$. But, when verifying (\ref{weak star}) for a fixed
$\gamma$ in this paper, 
a subset $\Omega_\gamma\subseteq{\mathscr N}_{div}^c$ is constructed
such that (\ref{weak star}) holds for all $Y\in\Omega_\gamma$ and 
$\frac{\partial}{\partial u}(Y^2)$
given by (\ref{first def}) on $\Omega_\gamma$ is the same 
for all even mollifiers $J$.

If $Y^\varepsilon$ approximates $Y$ then the standard method
for showing that $Y$ satisfies (\ref{weak star}) with
$\frac{\partial}{\partial u}(Y_t^2)$ defined by (\ref{first def})
would be: 
\begin{equation}\label{weak control}
\mbox{control\,\,
$\ee_\varepsilon\hspace{-2pt}\left[\rule{0pt}{13pt}\right.\hspace{-2pt}
\int_0^T\hspace{-5pt}\int_\NR\,\frac{\partial}{\partial u}\phi(s,u)\,
(Y_s^\varepsilon\star J_N)^2(u)\,\dd u\dd s
\hspace{-2pt}\left.\rule{0pt}{13pt}\right]^2$\,\,
in $\varepsilon,N,\phi$.}
\end{equation}
A good control of this type has been obtained in \cite{JG2010}
for $\phi={\bf 1}_{[r,t]}\otimes G$
using the density fluctuations
$Y^\varepsilon$ in $\sqrt{\varepsilon}$-asymmetric
exclusion as approximation scheme. But, sharp bounds
on the spectral gap of the symmetric exclusion processes
restricted to finite boxes were required.

In this paper it is demonstrated that (\ref{weak control}) can be
based on the weaker estimates obtained in \cite[Lemma 3.3]{A2012}. 
Using these weaker estimates makes it more difficult to verify 
that $Y$ satisfies (\ref{weak star}).
But the proof of Proposition \ref{cont_mart}, which is the main
achievement of this paper, presents a method of how to overcome this
difficulty.
Having a method based on weaker estimates might be beneficial when
it comes to a similar problem with other highly singular SPDEs.

Finally
it should be mentioned that the estimates used in this paper,
just as the estimates found in \cite{JG2010},
are only justified in the case where $Y$ is the
spatial derivative of the Cole-Hopf solution 
starting from Gaussian white noise on $\NR$
which is a stationary state.
In this case, in particular since this invariant state is Gaussian,
the state space of $Y$ can be relaxed to be
$D([0,T];{\mathscr S}^\prime(\NR))$
with ${\mathscr S}^\prime(\NR)$ being the space of tempered
distributions---see Remark \ref{prop measure}(i).
But, in the non-stationary case, the growth conditions implied by the
theorems in \cite{BG1997} would not allow for ${\mathscr S}^\prime(\NR)$
without further analysis. As a consequence 
$\frac{\partial}{\partial u}(Y^2)$ was defined to be an element of 
${\mathscr D}^\prime((0,T)\times\NR)$ 
to leave room for non-stationary solutions.

It remains an open problem to show that the spatial derivative of the
Cole-Hopf solution starting from initial conditions 
other than Gaussian white noise on $\NR$
satisfies (\ref{weak star}).\\

\noindent
{\it Acknowledgement}. 
The author thanks Martin Hairer for valuable comments.
\section{Notation and Results}\label{notations}
The approximation scheme for the conservative KPZ equation used in
this paper goes back to \cite{BG1997}. It is based on 
$\sqrt{\varepsilon}$-asymmetric exclusion processes and will be
briefly explained in what follows. The reader is referred
to \cite{L1999} for the underlying theory of exclusion processes.

Fix $\gamma\not=0$ and consider a scaling parameter $\varepsilon>0$
small enough such that $\sqrt{\varepsilon}\gamma\in[-1,1]$.
Denote by $(\Omega,{\cal F},\pp^\varepsilon_{\!\!\eta},\eta\in\{0,1\}^\NZ,
(\eta_t)_{t\ge 0})$ the strong Markov Feller process
whose generator $L_\varepsilon$ 
acts on local functions $f:\{0,1\}^\NZ\to\NR$ as
\begin{equation}\label{operator}
\begin{array}{rcl}
\displaystyle
\hspace{-2cm}
L_\varepsilon f(\eta)&=&\displaystyle
\sum_{x\in\NZ}\left(\rule{0pt}{12pt}\right.
(1+\sqrt{\varepsilon}\gamma)\,\eta(x)(1-\eta(x+1))[f(\eta^{x,x+1})-f(\eta)]\\
&&\hspace{3cm}\displaystyle
+\;(1-\sqrt{\varepsilon}\gamma)\,\eta(x)(1-\eta(x-1))[f(\eta^{x,x-1})-f(\eta)]
\left.\rule{0pt}{12pt}\right)
\end{array}
\end{equation}
where $\eta^{x,y}$ is standard notation for the operation
which exchanges the `spins' at $x$ and $y$.

Denote by $\nu_{1/2}$ the Bernoulli product  \label{initial}
measure on $\{0,1\}^{\NZ}$ satisfying
$\nu_{1/2}(\eta(x)=1)=1/2$ for all $x\in\NZ$. Define
$$\ppe\,=\,\int\pp^\varepsilon_{\!\!\eta}\,\dd\nu_{1/2}(\eta),
\quad
\xi_t(x)\,=\,\frac{\eta_t(x)-1/2}{\sqrt{1/4}},\;t\ge 0,\,x\in\NZ,$$
and notice that
the process $(\xi_t)_{t\ge 0}$ is a mean-zero stationary process on 
$(\Omega,{\cal F},\ppe)$ which takes values in $\{-1,1\}^{\NZ}$.

Denote by $\delta_{\varepsilon x}$ the Dirac measure
concentrated in the macroscopic point $\varepsilon x$
and define by
$$Y_t^\varepsilon\,=\,\sqrt{\varepsilon}
\sum_{x\in\NZ}\xi_{t\varepsilon^{-2}}(x)
\delta_{\varepsilon x},\quad t\ge 0,$$
the measure-valued density fluctuation field.
Fix a finite time horizon $T$ and regard
$Y^\varepsilon=(Y_t^\varepsilon)_{t\in[0,T]}$ as a random variable
taking values in the space $D([0,T];{\mathscr S}^\prime(\NR))$
of all cadlag functions mapping $[0,T]$ 
into the space of tempered distributions ${\mathscr S}^\prime(\NR)$.
Equip $D([0,T];{\mathscr S}^\prime(\NR))$ 
with the Skorokhod topology $J_1$ and let 
$Y$ be the notation for both
an element in and the identity map on $D([0,T];{\mathscr S}^\prime(\NR))$.
So $Y=(Y_t)_{t\in[0,T]}$ plays the role of the coordinate process on 
$D([0,T];{\mathscr S}^\prime(\NR))$ 
and it is evident that the topological $\sigma$-algebra on
$D([0,T];{\mathscr S}^\prime(\NR))$
is equal to 
${\cal F}^Y_T=\sigma(\{Y_t(G):t\le T,G\in{\mathscr S}(\NR)\})$.
\begin{theo}\label{bigBG}
{\rm \cite[Th.B.1 \& Prop.B.2]{BG1997}}
Let $\hat{\pp}_{\hspace{-3pt}\varepsilon}$ denote
the push forward of $\ppe$ with respect to the map $Y^\varepsilon$.
Then, when $\varepsilon\downarrow 0$,
the probability measures
$\hat{\pp}_{\hspace{-3pt}\varepsilon}$
converge weakly to a probability measure on $D([0,T];{\mathscr S}^\prime(\NR))$ 
which is denoted by $\ppo$ in what follows.
The measure $\ppo$ has the following properties:
\begin{itemize}
\item[(i)]
the support of $\ppo$ is a subset of $C([0,T];{\mathscr S}^\prime(\NR))$; 
\item[(ii)]
the process $Y$ is stationary under $\ppo$ satisfying
$Y_t\sim\mu,\,t\in[0,T]$, where $\mu$ is the mean-zero
Gaussian white noise measure with covariance
$\eeo Y_t(G)Y_t(H)=\int_\NR GH\,\dd u$;
\item[(iii)]
$\ppo$ is equal to the law of the spatial derivative of the so-called
Cole-Hopf solution of the KPZ equation for growing interfaces
starting from a two-sided Brownian motion.
\end{itemize}
\end{theo}
\begin{rem}\rm\label{prop measure}
\begin{itemize}
\item[(i)] The space used in Th.B.1 of \cite{BG1997} is
$D([0,T];{\mathscr D}^\prime(\NR))$. But this can be
relaxed to $D([0,T];{\mathscr S}^\prime(\NR))$ 
because $\nu_{1/2}$ is the initial condition of
$(\eta_t)_{t\ge 0}$.
Indeed, this implies that condition (2.13) 
on page 578 in \cite{BG1997} is satisfied
for $m\equiv 0$ and one can rule out that the functions $f_X$ used in
the proof of Th.B.1 have exponential growth.
 
\item[(ii)] This result in \cite{BG1997} is stronger than the tightness of 
$\{\hat{\pp}_{\hspace{-3pt}\varepsilon},\,\varepsilon>0\}$
shown in \cite{JG2010} as
tightness would only give the weak convergence 
with respect to certain subsequences
$\varepsilon_k,\,\varepsilon_k\downarrow 0$, with possibly different
limit points.
The identification of all limit points is a consequence of the
Cole-Hopf transform for discrete systems applied in \cite{BG1997}.
\end{itemize}
\end{rem}
\begin{defi}\rm
The coordinate process $Y$ on the probability space
$(D([0,T];{\mathscr S}^\prime(\NR)),{\cal F}_T^Y,\ppo)$
is called 
{\it Cole-Hopf solution of the conservative KPZ equation}
(\ref{star}).
\end{defi}

The following two results 
whose proofs will be given in the next section
form the basis for the method of verification
used in this paper to show
that $Y$ solves equation (\ref{star}) in the sense of (\ref{weak star})
where $\frac{\partial}{\partial u}(Y_t^2)$ 
is defined by (\ref{first def}).
Notice that, by technical reasons, 
the mollifier $J\in{\mathscr D}(\NR)$ 
defining $J_N$ by $u\mapsto NJ(Nu),\,N\ge 1$, 
should be taken to be \underline{even}.
\begin{lemma}\label{cauchy}
Fix $G\in{\mathscr S}(\NR)$. Then
\begin{eqnarray*}
&&\int_0^T\dd t\,\eeo
\left[
\int_{0}^{t}\hspace{-5pt}\int_\NR 
G^\prime(u)\,
\left(\rule{0pt}{12pt}
(Y_{s}\star J_{\tilde{N}})^2(u)
-
(Y_s\star J_{N})^2(u)\right)
\,\dd u\dd s\right]^2\\
&\le&\rule{0pt}{30pt}
e^T\CJ\,N^{-1/3}
\sum_{m=1}^3\sup_{u}
|(1+u^2)\frac{\partial^m}{\partial u^m}G(u)|^2
\end{eqnarray*}
for all $\tilde{N}\ge N\ge 1$
where $\CJ$ is a constant which only depends on the choice of
the mollifier $J$.
\end{lemma}
This lemma is proven using the estimates obtained in \cite[Lemma 3.3]{A2012} 
by applying a resolvent-type method. It only gives a bound on 
the $(\ell\otimes\ppo)$\,-\,average of the square of the functional
$$(t,Y)\mapsto
\int_{0}^{t}\hspace{-5pt}\int_\NR
G^\prime(u)\,
\left(\rule{0pt}{12pt}
(Y_{s}\star J_{\tilde{N}})^2(u)
-
(Y_s\star J_{N})^2(u)\right)
\dd u\dd s$$
where $\ell$ denotes the Lebesgue measure on $[0,T]$.
The main disadvantage of using an $L^2(\ell\otimes\ppo)$\,-\,estimate 
of the above  functional is that it complicates the method of 
identifying the Brownian sheet in (\ref{weak star}). 
The next proposition deals with each single step of this method in
detail. It's proof is also based on \cite[Lemma 3.3]{A2012}, only.
This means that fairly weak $L^2(\ell\otimes\pp)$ a priori estimates
are still good enough for solving singular SPDEs.

Define the map
$$\mathfrak{M}_N:
D([0,T];{\mathscr S}^\prime(\NR))\to D([0,T];{\mathscr S}^\prime(\NR))$$
by
$$\mathfrak{M}_N(Y)_t^G\,=\,
Y_t(G)-Y_0(G)
-\int_0^t\,Y_s(G^{\prime\prime})\,\dd s
\,+{\gamma}
\int_{0}^{t}\hspace{-5pt}\int_\NR
G^\prime(u)
(Y_s\star J_{N})^2(u)
\,\dd u\dd s.$$
Applying Lemma \ref{cauchy} 
gives that, for every $G\in{\mathscr S}(\NR)$, there exists a 
${\cal B}([0,T])\otimes{\cal F}^Y_T$\,-\,measurable process
$$\tilde{M}^G:
[0,T]\times D([0,T];{\mathscr S}^\prime(\NR))\to\NR$$
such that 
\begin{equation}\label{converge}
\int_0^T\dd t\,\eeo
\left[\rule{0pt}{12pt}\tilde{M}^G_t-\mathfrak{M}_N(Y)_t^G\right]^2
\,\to\,0,\quad N\to\infty.
\end{equation}
Denote by $\NF$ the filtration $({\cal F}_t)_{t\in[0,T]}$ with 
${\cal F}_t=
\sigma(\{Y_s(G):s\le t,G\in{\mathscr S}(\NR)\}\cup{\cal N})$ 
where ${\cal N}$ is the
collection of all $\ppo$-null sets in ${\cal F}^Y_T$.
\begin{prop}\label{cont_mart}
\begin{itemize}
\item[(i)] For every $G\in{\mathscr S}(\NR)$,
there exists an $\NF$-adapted process
$M^G=(M^G_t)_{t\in[0,T]}$
on $(D([0,T];{\mathscr S}^\prime(\NR)),{\cal F}^Y_T,\ppo)$
which is a continuous version of $\tilde{M}^G$ in the following sense:
there is a measurable subset ${\cal T}_G\subseteq[0,T]$ with
$\ell({\cal T}_G)=T$ such that $\tilde{M}_t^G=M^G_t$ a.s.\ for all
$t\in{\cal T}_G$.
For every positive $T^\prime<T$, when restricted to $[0,T^\prime]$,
the process $M^G$ is a square integrable $\NF$-\,martingale.
\item[(ii)] For every $G\in{\mathscr S}(\NR)$,
the process $M^G=(M^G_t)_{t\in[0,T]}$  is an $\NF$-Brownian motion
with variance $2\|G^\prime\|^2_2$ on the probability space
$(D([0,T];{\mathscr S}^\prime(\NR)),{\cal F}^Y_T,\ppo)$.

\item[(iii)] It holds that
$$M_t^{a_1 G_1+a_2 G_2}\,=\,a_1 M_t^{G_1}+a_2 M_t^{G_2}
\quad\mbox{a.s.}$$
for every $t\in[0,T],\,a_1,a_2\in\NR$ and $G_1,G_2\in{\mathscr S}(\NR)$.

\item[(iv)] The process $M_t^G$ 
indexed by $t\in[0,T]$ and $G\in{\mathscr S}(\NR)$
is a centred Gaussian process on 
$(D([0,T];{\mathscr S}^\prime(\NR)),{\cal F}^Y_T,\ppo)$
with covariance
$$\eeo M_{t_1}^{G_1}M_{t_2}^{G_2}\,=\,
2(t_1\wedge t_2)\int_\NR G_1^\prime(u)G_2^\prime(u)\,\dd u$$
hence there is a Brownian sheet $B(t,u),\,t\in[0,T],\,u\in\NR$, on
$(D([0,T];{\mathscr S}^\prime(\NR)),{\cal F}^Y_T,\ppo)$
such that
$$M_t^G\,=\,\sqrt{2}
\int_\NR B(t,u)G^{\prime\prime}(u)\,\dd u
\quad\mbox{a.s.}$$
for every $t\in[0,T]$ and $G\in{\mathscr S}(\NR)$.
\end{itemize}
\end{prop}

In what follows let $M=(M_t)_{t\in[0,T]}$ denote the continuous
${\mathscr S}^\prime(\NR)$\,-\,valued process defined by
\begin{equation}\label{M is sheet}
M_t(G)\,\stackrel{\mbox{\tiny def}}{=}\,\sqrt{2}
\int_\NR B(t,u)G^{\prime\prime}(u)\,\dd u,
\quad t\in[0,T],\,G\in{\mathscr S}(\NR).
\end{equation}
Remark that, by Schwartz' kernel theorem,
$M$ and $Y$ can also be considered random variables taking values in 
${\mathscr D}^\prime((0,T)\times\NR)$ such that
$$\int_0^T\hspace{-5pt}\dd t\,g^\prime(t)
\!\left[
-Y_t(G)+Y_0(G)
+\int_0^t\hspace{-3pt}Y_s(G^{\prime\prime})\,\dd s
\,+M_t(G)
\right]
\,=\,
\langle g\otimes G\,,
\frac{\partial}{\partial t}Y-\frac{\partial^2}{\partial u^2}Y
-\frac{\partial}{\partial t}M\rangle$$
for all $g\in{\mathscr D}((0,T)),\,G\in{\mathscr D}(\NR)$
where $\langle\cdot\,,\cdot\rangle$ denotes the dual pairing between
${\mathscr D}((0,T)\times\NR)$ and ${\mathscr D}^\prime((0,T)\times\NR)$.
Notice that the last equality can be extended to hold for all
$g\in C^1([0,T])$ with $g(T)=0$ and $G\in{\mathscr S}(\NR)$.
Then it is an easy consequence of Lemma \ref{cauchy}, (\ref{converge}) and the
Cauchy-Schwarz inequality that
\begin{eqnarray}\label{not energy inequality}
&&\eeo\left|
-\!\int_0^T\hspace{-7pt}\int_\NR
g(t)G^\prime(u)\,(Y_t\star J_{N})^2(u)\,\dd u\dd t
-\langle g\otimes G\,,
\frac{\partial}{\partial t}Y-\frac{\partial^2}{\partial u^2}Y
-\frac{\partial}{\partial t}M\rangle/{\gamma}
\right|^2\nonumber\\
&\le&\rule{0pt}{18pt}
2e^T\CJ\,N^{-1/3}\,
\|g^\prime\|_{L^2[0,T]}^2
\sum_{m=1}^3\sup_{u}
|(1+u^2)\frac{\partial^m}{\partial u^m}G(u)|^2
,\quad N\ge 1,
\end{eqnarray}
for all $g\in C^1([0,T])$ with $g(T)=0$ and $G\in{\mathscr S}(\NR)$.

The next step consists in finding a subsequence $(N_k)_{k=1}^\infty$
and a subset $\Omega_\gamma\in{\cal F}_T^Y$ of measure
$\ppo(\Omega_\gamma)=1$ such that 
$\Omega_\gamma\subseteq{\mathscr N}_{div}^c$
where ${\mathscr N}_{div}$ is defined by (\ref{Ndiv}).
The ultimate goal would of course be a subsequence
$(N_k)_{k=1}^\infty$ which is the same for all $\gamma\not=0$.

For this purpose it turns out to be useful to think of the function
$$(t,u)\mapsto(Y_t\star J_N)^2(u)
\quad\mbox{where}\quad
Y\in D([0,T];{\mathscr S}^\prime(\NR))$$
as a regular distribution in 
${\mathscr D}^\prime((0,T)\times\NR)$. 
This regular distribution is denoted  by $(Y\star_2 J_N)^2$
in what follows. Notice the notation $\star_2$ which emphasises that the
convolution only acts on the space component of $Y$.

Then the idea is to construct a 
Banach space $({\cal E},|\!|\!|\cdot|\!|\!|)$
satisfying
${\mathscr D}((0,T)\times\NR)\subseteq{\cal E}^\prime\subseteq
L^2([0,T]\times\NR)\subseteq{\cal E}$
such that
\begin{equation}\label{main esti}
\eeo|\!|\!|\frac{\partial}{\partial u}(Y\star_2 J_N)^2- 
(\frac{\partial}{\partial t}Y-\frac{\partial^2}{\partial u^2}Y
-\frac{\partial}{\partial t}M)/{\gamma}|\!|\!|^2
\,\le\,const/N^\alpha,
\quad N\ge 1,
\end{equation}
for some $\alpha>0$.
\begin{rem}\rm\label{proof theo}
\begin{itemize}\item[]\hspace{-25pt}
Suppose for now that (\ref{main esti}) can be achieved
by finding $|\!|\!|\cdot|\!|\!|,\alpha$ and $const$
where the latter might depend on $T,J$ and $\gamma$.
Choosing $(N_k)_{k=1}^\infty$ to be
$$N_k\,=\left\{\begin{array}{lcc}
\mbox{$k^{\tilde{\alpha}}$\, for some
$\tilde{\alpha}>1/\alpha$}&:&\alpha\le 1\\
k&:&\alpha>1\rule{0pt}{13pt}
\end{array}\right.$$
would then yield
$$\sum_{k=1}^\infty
\ppo(\{|\!|\!|\frac{\partial}{\partial u}(Y\star_2 J_{N_k})^2-
(\frac{\partial}{\partial t}Y-\frac{\partial^2}{\partial u^2}Y
-\frac{\partial}{\partial t}M)/{\gamma}|\!|\!|
\,\ge\,\delta\})\,<\,\infty,
\quad\forall\,\delta>0,$$
hence
$$|\!|\!|\frac{\partial}{\partial u}(Y\star_2 J_{N_k})^2-
(\frac{\partial}{\partial t}Y-\frac{\partial^2}{\partial u^2}Y
-\frac{\partial}{\partial t}M)/{\gamma}|\!|\!|
\longrightarrow 0,
\quad k\to\infty,$$
for all $Y\in\Omega_{{\gamma}}$ for some 
$\Omega_{{\gamma}}\in{\cal F}^Y_T$ with 
$\ppo(\Omega_{{\gamma}})=1$.
Since weak convergence in ${\cal E}$ implies
weak convergence in ${\mathscr D}^\prime((0,T)\times\NR)$
one would obtain that 
\begin{eqnarray*}
\langle\phi\,,
\frac{\partial}{\partial t}Y-\frac{\partial^2}{\partial u^2}Y
-\frac{\partial}{\partial t}M\rangle/{\gamma}
&=&
\lim_{k\to\infty}\langle\phi\,,
\frac{\partial}{\partial u}(Y\star_2 J_{N_k})^2\rangle\\
&=&\rule{0pt}{20pt}
-\lim_{k\to\infty}\int_0^T\hspace{-7pt}\int_\NR
\frac{\partial}{\partial u}\phi(t,u)
\,(Y_t\star J_{N_k})^2(u)\,\dd u\dd t
\end{eqnarray*}
for all $\phi\in{\mathscr D}((0,T)\times\NR)$ and
$Y\in\Omega_{{\gamma}}$ which obviously means
$\Omega_\gamma\subseteq{\mathscr N}_{div}^c$
where the chosen subsequence $(N_k)_{k=1}^\infty$ 
would indeed be the same for all $\gamma\not=0$.
Notice that 
$\langle\phi\,,
\frac{\partial}{\partial t}Y-\frac{\partial^2}{\partial u^2}Y
-\frac{\partial}{\partial t}M\rangle/{\gamma}$
does not depend on the choice of $J$ so that
$\frac{\partial}{\partial u}(Y^2)$
given by (\ref{first def}) on $\Omega_\gamma$ would be the same 
for all even mollifiers $J$.
Furthermore, the equality in (\ref{weak star})
would also be true for all 
$\phi\in{\mathscr D}((0,T)\times\NR)$ and all $Y\in\Omega_\gamma$
because 
$\partial{M}/\partial{t}=
\sqrt{2}\,{\partial^3 B}/\partial{t}/{\partial u^2}$
by (\ref{M is sheet}).
\end{itemize}
\end{rem}

So it remains to justify (\ref{main esti}).
Of course, one wants to use the bounds 
given by the right-hand side of (\ref{not energy inequality}) 
to construct the Banach space $({\cal E},|\!|\!|\cdot|\!|\!|)$
but some care is needed to ensure that 
${\mathscr D}((0,T)\times\NR)\subseteq{\cal E}^\prime$.
A straight forward approach to tackle this problem is using a
so-called negative-order Sobolev space which is introduced next.

First observe that
\begin{equation}\label{sup by L2}
\sup_{u\in\NR}|(1+u^2)H(u)|^2
\,\le\,
4\,\|(1+u^2)H\|^2_{L^2(\NR)}+
2\,\|(1+u^2)H\|_{L^2(\NR)}\,\|(1+u^2)H^\prime\|_{L^2(\NR)}
\end{equation}
for any test function $H\in{\mathscr S}(\NR)$.
Now let $(g_m)_{m=1}^\infty$ be the eigenbasis 
of the one-dimensional Laplacian on $[0,T]$ with Dirichlet boundary
conditions and let $(G_n)_{n=1}^\infty$ be the
collection of Hermite functions.
Then $(g_n\otimes G_m)_{n,m}$ forms an orthonormal basis in
$L^2([0,T]\times\NR)$ and it follows from 
(\ref{not energy inequality}) and (\ref{sup by L2}) that
\begin{equation}\label{for space}
\eeo\left|\rule{0pt}{12pt}\right.
\langle g_m\otimes G_n\,,
\frac{\partial}{\partial u}(Y\star_2 J_{N_k})^2- 
(\frac{\partial}{\partial t}Y-\frac{\partial^2}{\partial u^2}Y
-\frac{\partial}{\partial t}M)/{\gamma}\rangle
\left.\rule{0pt}{12pt}\right|^2
\,\le\,
{const}\cdot m^2 n^6/N^{1/3}
\end{equation}
where $const$ does not depend on the choice of $m$ and $n$.
Of course the factor $m^2$ goes back to the
eigenvalue associated with $g_m$ and,
using the combinatorical properties of the Hermite functions,
$O(n^6)$ is a quite crude estimate
of the norms of $H$ and its derivative in (\ref{sup by L2})
when $H=G_n^\prime,G_n^{\prime\prime},G_n^{\prime\prime\prime}$.
So an appropriate choice of the Banach space ${\cal E}$ 
is the completion of ${\mathscr D}((0,T)\times\NR)$ 
with respect to the norm $|\!|\!|\cdot|\!|\!|$ given by
$$|\!|\!|\phi|\!|\!|^2\,=\,
\sum_{m,n}\left[(m^3+n^3)m^2 n^6\right]^{-1}\,
\langle g_m\otimes G_n\,,\phi\rangle^2.$$
Notice that ${\mathscr D}((0,T)\times\NR)\subseteq{\cal E}^\prime$
is a standard consequence when choosing
$(g_m)_{m=1}^\infty$ and $(G_n)_{n=1}^\infty$ as above.

Using this Banach space and applying (\ref{for space}) to calculate
$\eeo|\!|\!|\frac{\partial}{\partial u}(Y\star_2 J_N)^2- 
(\frac{\partial}{\partial t}Y-\frac{\partial^2}{\partial u^2}Y
-\frac{\partial}{\partial t}M)/{\gamma}|\!|\!|^2$ 
results in (\ref{main esti}) for $\alpha=1/3$
hence Remark \ref{proof theo} proves the following theorem.
\begin{theo}\label{pde}
\begin{itemize}\item[(i)]
There exists a subsequence $(N_k)_{k=1}^\infty$ such that
for every ${\gamma}\not=0$ there is a set
$\Omega_{{\gamma}}\in{\cal F}^Y_T$ with 
$\ppo(\Omega_{{\gamma}})=1$ such that
$\Omega_\gamma\subseteq{\mathscr N}_{div}^c$
where ${\mathscr N}_{div}$ is defined by (\ref{Ndiv})
and $\frac{\partial}{\partial u}(Y^2)$
given by (\ref{first def}) on $\Omega_\gamma$ 
is the same for all even mollifiers $J$.
\item[(ii)]
There exists a Brownian sheet $B(t,u),\,t\in[0,T],\,u\in\NR$, 
on 
$(D([0,T];{\mathscr D}^\prime(\NR)),$ ${\cal F}^Y_T,\ppo)$ such that
the coordinate process $Y$ solves the equation (\ref{star}) in the sense
of (\ref{weak star}).
\end{itemize}
\end{theo}
\begin{rem}\rm
\begin{itemize}\item[(i)]
The choice of the subsequence used in the definition 
(\ref{first def}) of $\frac{\partial}{\partial u}(Y^2)$
depends on the power $\alpha$ needed to establish (\ref{main esti}).
The power $\alpha=1/3$ used in this paper goes back to 
\cite{A2007}. The results in \cite{JG2010} suggest that $\alpha=1/2$
seems to be possible. However, for the purpose of giving rigorous
sense to the equation (\ref{star}), the choice of an optimal
subsequence is not intrinsic and so the author used what he had proved
himself in \cite{A2007}. 
But, in the light of the new techniques applied in \cite{H2012},
he would like to \underline{conjecture} the following: 
equation (\ref{star}) holds true in the sense of (\ref{weak star})
using $N_k=k$ in the definition 
(\ref{first def}) of $\frac{\partial}{\partial u}(Y^2)$.
\item[(ii)]
It is a consequence of Theorem \ref{pde}(i) that
$\bigcup_{{\gamma}\not=0}\Omega_{{\gamma}}
\subseteq{\mathscr N}_{div}^c$.
But, as shown in \cite{BG1997}, each measure $\ppo$
is related to the solution of 
a corresponding stochastic heat equation
$$\frac{\partial}{\partial t}Z
\,=\,
\frac{\partial^2}{\partial u^2}Z
+\gamma\sqrt{2}\,Z\,\frac{\partial^2 B}{\partial t\partial u},
\quad\gamma\not=0,$$
through the Cole-Hopf transform and changing the diffusion coefficient 
by $\gamma$ indicates that all measures
$\ppo,\,\gamma\not=0$, are singular to each other.
Thus, the set $\bigcup_{{\gamma}\not=0}\Omega_{{\gamma}}$
is not too small since
$\ppo(\Omega_{{\gamma}})=1$ for all $\gamma\not=0$.
\end{itemize}
\end{rem}
\section{Proofs}
This section contains the proofs of Lemma \ref{cauchy}
and Proposition \ref{cont_mart}, but first, further notation and
auxiliary results need to be provided.

Fix $\varepsilon>0$ small enough such that
$\sqrt{\varepsilon}\gamma\in[-1,1]$,
fix a test function $G\in{\mathscr S}(\NR)$
and denote by $\|\cdot\|_p$ the norm in $L^p(\NR),\,1\le p\le\infty$.
Then $$M^{G,\varepsilon}_t\,=\,Y_t^\varepsilon(G)-Y_0^\varepsilon(G)-
\int_0^t \varepsilon^{-2}L_\varepsilon Y_s^\varepsilon(G)
\,\dd s,\;t\ge 0,$$
is a martingale on $(\Omega,{\cal F},\ppe)$
by standard theory on strong Markov processes and
\begin{equation}\label{gen on field}
\int_0^t \varepsilon^{-2}L_\varepsilon Y_s^\varepsilon(G)\,\dd s
\,=\int_0^t\varepsilon^{-\frac{3}{2}}
\sum_{x\in\NZ}G(\varepsilon x)
L_\varepsilon\xi_{s\varepsilon^{-2}}(x)\,\dd s,\;t\ge 0,
\end{equation}
where
\begin{equation}\label{generator}
\begin{array}{rcl}
\displaystyle
L_\varepsilon\xi_{s\varepsilon^{-2}}(x)
&=&\left[\rule{0pt}{12pt}
\left(\rule{0pt}{10pt}
\xi_{s\varepsilon^{-2}}(x-1)-2\xi_{s\varepsilon^{-2}}(x)
+\xi_{s\varepsilon^{-2}}(x+1)
\right)\right.\\
&&\hspace{17pt}\displaystyle
\left.+\,\sqrt{\varepsilon}\gamma\left(\rule{0pt}{10pt}
\xi_{s\varepsilon^{-2}}(x)\xi_{s\varepsilon^{-2}}(x+1)
-\xi_{s\varepsilon^{-2}}(x-1)\xi_{s\varepsilon^{-2}}(x)\right)
\rule{0pt}{12pt}\right]
\end{array}
\end{equation}
follows from (\ref{operator}).
Substituting (\ref{generator}) into (\ref{gen on field}),
performing a summation by parts 
and approximating by Taylor expansion implies
\begin{eqnarray*}
\int_0^t \varepsilon^{-2}L_\varepsilon Y_s^\varepsilon(G)\,\dd s
&=&
\int_0^t Y_s^\varepsilon(G^{\prime\prime})\,\dd s
\,-\,\gamma\int_0^t
\sum_{x\in\NZ}G^\prime(\varepsilon x)\,
\xi_{s\varepsilon^{-2}}(x)\xi_{s\varepsilon^{-2}}(x+1)\,\dd s\\
&+&\frac{\gamma\varepsilon}{2}\int_0^t
\sum_{x\in\NZ}G^{\prime\prime}(\varepsilon x)\,
\xi_{s\varepsilon^{-2}}(x)\xi_{s\varepsilon^{-2}}(x+1)\,\dd s
\,+R^G_\varepsilon(t)
\end{eqnarray*}
with
\begin{equation}\label{taylor esti}
|R^G_\varepsilon(t)|\,\le\,
\sqrt{\varepsilon}\,\frac{1}{6}
(2+\sqrt{\varepsilon}\gamma)
(\pi+2\varepsilon)
\|(1+u^2)G^{\prime\prime\prime}\|_\infty\cdot t,
\quad t\ge 0,
\end{equation}
where $\pi+2\varepsilon$ is an upper bound of the discretization
of the integral $\int_\NR(1+u^2)^{-1}\dd u$ in this context.
Now, by notational purpose, define
\begin{equation}\label{RG0}
R_{\varepsilon,N}^{G^\prime,0}(t)
\,\stackrel{\mbox{\tiny def}}{=}\,
\frac{\varepsilon}{2}\int_0^{t}
\sum_{x\in\NZ}G^{\prime\prime}(\varepsilon x)
\xi_{s\varepsilon^{-2}}(x)\xi_{s\varepsilon^{-2}}(x+1)\,\dd s,
\quad t\ge 0,
\end{equation}
although the right-hand side does not depend on $N$ 
and includes $G^{\prime\prime}$ instead of $G^\prime$.
Using this notation leads to the decomposition
\begin{equation}\label{approxi}
\begin{array}{rcl}
\displaystyle
M^{G,\varepsilon}_t+\,R^G_\varepsilon(t)
\,+\,\gamma R_{\varepsilon,N}^{G^\prime,0}(t)
&=&\displaystyle
Y_t^\varepsilon(G)-Y_0^\varepsilon(G)
-\int_0^t Y_s^\varepsilon(G^{\prime\prime})\,\dd s\\
&+&\displaystyle
{\gamma}\int_0^t
\sum_{x\in\NZ}G^{\prime}(\varepsilon x)
\xi_{s\varepsilon^{-2}}(x)\xi_{s\varepsilon^{-2}}(x+1)\,\dd s
\end{array}
\end{equation}
for all $t\ge 0$.

It turns out to be useful 
to rewrite the difference below as follows
\begin{equation}\label{rewrite diff}
\int_{0}^{t}\hspace{-5pt}\int_\NR
G^\prime(u)
(Y_s^\varepsilon\star J_{N})^2(u)
\,\dd u\dd s
\,-\,
\int_0^t
\sum_{x\in\NZ}G^\prime(\varepsilon x)
\xi_{s\varepsilon^{-2}}(x)\xi_{s\varepsilon^{-2}}(x+1)\,\dd s
\,=
\sum_{i=1}^4 R_{\varepsilon,N}^{G^\prime,i}(t)
\end{equation}
where
\begin{equation}\label{RGi}
R_{\varepsilon,N}^{G^\prime,i}(t)
\,\stackrel{\mbox{\tiny def}}{=}
\int_0^{t}
V_{\varepsilon,N}^{G^\prime,i}(\xi_{s\varepsilon^{-2}})\,\dd s,
\quad t\ge 0,\quad i=1,2,3,4,
\end{equation}
are given by
\begin{eqnarray*}
V_{\varepsilon,N}^{G^\prime,1}(\xi)&=&
\sum_{x\in\NZ}
\int_\NR
[G^\prime(u)-G^\prime(\varepsilon x)]
J_N(u-\varepsilon x)
\sum_{\tilde{x}\in\NZ}
\varepsilon
J_N(u-\varepsilon\tilde{x})\,\dd u\;
\xi(x)\xi(\tilde{x}),\\
V_{\varepsilon,N}^{G^\prime,2}(\xi)&=&
\varepsilon
\sum_{x\in\NZ}G^\prime(\varepsilon x)
\int_\NR J_N^2(u-\varepsilon x)\,\dd u\;
\xi(x)[\xi(x)-\xi(x+1)],\\
V_{\varepsilon,N}^{G^\prime,3}(\xi)&=&
\varepsilon
\sum_{x\not=\tilde{x}}G^\prime(\varepsilon x)
\int_\NR J_N(u-\varepsilon x)J_N(u-\varepsilon\tilde{x})\,\dd u\;
\xi(x)[\xi(\tilde{x})-\xi(x+1)],\\
V_{\varepsilon,N}^{G^\prime,4}(\xi)&=&
\sum_{x\in\NZ}G^\prime(\varepsilon x)
\int_\NR J_N(u-\varepsilon x)
\hspace{-3pt}\left[
\sum_{\tilde{x}\in\NZ}\varepsilon J_N(u-\varepsilon\tilde{x})-1
\right]\!\dd u\;
\xi(x)\xi(x+1).
\end{eqnarray*}
Notice that $\int_\NR G^\prime(u)\dd u=0$ hence
the following lemma can be applied in this context.
\begin{lemma}\label{key result}
{\rm\cite[Lemma 3.3]{A2012}}
Recall (\ref{RGi}) for the definition of
$R_{\varepsilon,N}^{G^\prime,i},\,i=1,2,3,4$.
Then
$$\begin{array}{crcl}
\hspace{-2.5cm}(i)\hspace{.2cm}&\displaystyle
\int_0^T\hspace{-3pt}\dd t\,\ee_\varepsilon
\left[R_{\varepsilon,N}^{G^\prime,1}(t)\right]^2
&\le&\displaystyle
e^T\tilde{\CJ}
\left(\rule{0pt}{11pt}
\frac{\|(1+u^2)G^{\prime\prime\prime}\|_\infty^2}{N^2}
+\frac{\|(1+u^2)G^{\prime\prime}\|_\infty^2}{N}
\right)\\
\rule{0pt}{30pt}
\hspace{-2.5cm}(ii)\hspace{.2cm}&\displaystyle
\int_0^T\hspace{-3pt}\dd t\,\ee_\varepsilon
\left[R_{\varepsilon,N}^{G^\prime,2}(t)\right]^2
&\le&\displaystyle
e^T\tilde{\CJ}
\left(\rule{0pt}{11pt}
\varepsilon^2 N^2\,\|G^{\prime\prime}\|_\infty^2
+\varepsilon N^2\,\|(1+u^2)G^\prime\|_\infty^2
\right)\\
\rule{0pt}{30pt}
\hspace{-2.5cm}(iii)\hspace{.2cm}&\displaystyle
\int_0^T\hspace{-3pt}\dd t\,\ee_\varepsilon
\left[R_{\varepsilon,N}^{G^\prime,3}(t)\right]^2
&\le&\displaystyle
e^T\tilde{\CJ}
\,\frac{\|(1+u^2)G^\prime\|_\infty^2}{N^{1/3}}\\
\rule{0pt}{30pt}
\hspace{-2.5cm}(iv)\hspace{.2cm}&\displaystyle
\int_0^T\hspace{-3pt}\dd t\,\ee_\varepsilon
\left[R_{\varepsilon,N}^{G^\prime,4}(t)\right]^2
&\le&\displaystyle
T^3\tilde{\CJ}
\,\varepsilon^2 N^4\|G^\prime\|_1^2
\end{array}$$
for all $\varepsilon>0,N\ge 1$ where $\tilde{\CJ}$ is a constant
which only depends on the choice of the mollifier $J$.
\end{lemma}
\begin{rem}\rm
Recall the definition of $R_{\varepsilon,N}^{G^\prime,0}$
given in (\ref{RG0}) which does not depend on $N$.
Then the rate of convergence
\begin{equation}\label{extra bound}
\int_0^T\dd t\,\ee_\varepsilon\hspace{-3pt}
\left[R_{\varepsilon,N}^{G^\prime,0}(t)\right]^2
\,=\,O(\varepsilon^2),
\quad\mbox{$\varepsilon\downarrow 0$, uniformly in $N$,}
\end{equation}
follows from Remark 1(iii) in \cite{A2007}
by the same method 
used in the proof of the above lemma in \cite{A2012}.
\end{rem}
{\bf Proof of Lemma \ref{cauchy}.}
Fix $N\ge 1$, fix $\delta>0$ and choose
$N_\delta\ge N$ such that
$$8e^T\tilde{\CJ}\left(\rule{0pt}{12pt}\right.
N_\delta^{-2}\|(1+u^2)G^{\prime\prime\prime}\|_\infty^2
+N_\delta^{-1}\|(1+u^2)G^{\prime\prime}\|_\infty^2
+N_\delta^{-1/3}\|(1+u^2)G^\prime\|_\infty^2
\left.\rule{0pt}{12pt}\right)
\,\le\,\delta/4$$
where $\tilde{\CJ}$ is the constant 
appearing in Lemma \ref{key result}.
Then
\begin{eqnarray*}
&&\int_0^T\dd t\,\eeo
\left[
\int_{0}^{t}\hspace{-5pt}\int_\NR G^\prime(u)\,
\left(\rule{0pt}{12pt}
(Y_{s}\star J_{{N_\delta}})^2(u)
-
(Y_s\star J_{N})^2(u)\right)
\,\dd u\dd s\right]^2\\
&=&\int_0^T\dd t\,
\int_{0}^{t}\hspace{-5pt}\int_\NR
\int_{0}^{t}\hspace{-5pt}\int_\NR
\dd s_1\dd u_1\dd s_2\dd u_2\;
G^\prime(u_1)\,G^\prime(u_2)\\
&&\times\;\eeo\!
\left(\rule{0pt}{12pt}
(Y_{s_1}\star J_{{N_\delta}})^2(u_1)
-
(Y_{s_1}\star J_{N})^2(u_1)\right)
\left(\rule{0pt}{12pt}
(Y_{s_2}\star J_{{N_\delta}})^2(u_2)
-
(Y_{s_2}\star J_{N})^2(u_2)\right)
\end{eqnarray*}
where by Lemma \ref{uni bound} in the Appendix
$$\eeo\!
\left(\rule{0pt}{12pt}
(Y_{s_1}\star J_{{N_\delta}})^2(u_1)
-
(Y_{s_1}\star J_{N})^2(u_1)\right)
\left(\rule{0pt}{12pt}
(Y_{s_2}\star J_{{N_\delta}})^2(u_2)
-
(Y_{s_2}\star J_{N})^2(u_2)\right)$$
$$=\,\lim_{\varepsilon\downarrow 0}\,
\hat{\ee}_{\varepsilon}\!
\left(\rule{0pt}{12pt}
(Y_{s_1}\star J_{{N_\delta}})^2(u_1)
-
(Y_{s_1}\star J_{N})^2(u_1)\right)
\left(\rule{0pt}{12pt}
(Y_{s_2}\star J_{{N_\delta}})^2(u_2)
-
(Y_{s_2}\star J_{N})^2(u_2)\right)$$
such that
\begin{eqnarray*}
&&\left|\hat{\ee}_{\varepsilon}\!
\left(\rule{0pt}{12pt}
(Y_{s_1}\star J_{{N_\delta}})^2(u_1)
-
(Y_{s_1}\star J_{N})^2(u_1)\right)
\left(\rule{0pt}{12pt}
(Y_{s_2}\star J_{{N_\delta}})^2(u_2)
-
(Y_{s_2}\star J_{N})^2(u_2)\right)\right|^2\\
&\le&
\hat{f}(\|J_N\|^2,\|J_{{N_\delta}}\|^2)
\end{eqnarray*}
for all $\varepsilon\le 1,\,0\le s_1,s_2\le T$ and $u_1,u_2\in\NR$.
Hence, by dominated convergence, it follows that
\begin{eqnarray}\label{epsilon of N and tilde N}
&&\int_0^T\dd t\,\eeo
\left[
\int_{0}^{t}\hspace{-5pt}\int_\NR G^\prime(u)\,
\left(\rule{0pt}{12pt}
(Y_{s}\star J_{{N_\delta}})^2(u)
-
(Y_s\star J_{N})^2(u)\right)
\,\dd u\dd s\right]^2
\nonumber\\
&\le&
\frac{\delta}{2}\,+
\int_0^T\dd t\,\hat{\ee}_{\varepsilon}
\left[
\int_{0}^{t}\hspace{-5pt}\int_\NR G^\prime(u)\,
\left(\rule{0pt}{12pt}
(Y_{s}\star J_{{N_\delta}})^2(u)
-
(Y_s\star J_{N})^2(u)\right)
\,\dd u\dd s\right]^2
\end{eqnarray}
if $\varepsilon=\varepsilon_{N,{N_\delta}}>0$ 
is chosen to be sufficiently small. 

Using (\ref{rewrite diff}),
the last summand can be further estimated by
$$8\sum_{i=1}^4
\int_0^T\hspace{-5pt}\dd t\,
\ee_\varepsilon\left[
R_{\varepsilon,N}^{G^\prime,i}(t)
\right]^2
\,+\,
8\sum_{i=1}^4
\int_0^T\hspace{-5pt}\dd t\,
\ee_\varepsilon\left[
R_{\varepsilon,N_\delta}^{G^\prime,i}(t)
\right]^2$$
where
\begin{eqnarray*}
\sum_{i=1}^4
\int_0^\infty\hspace{-5pt}\dd t\,
\ee_\varepsilon\left[
R_{\varepsilon,N}^{G^\prime,i}(t)
\right]^2
&\le&
e^T\tilde{\CJ}\left(\rule{0pt}{12pt}\right.
\frac{\|(1+u^2)G^{\prime\prime\prime}\|_\infty^2}{N^2}
+\frac{\|(1+u^2)G^{\prime\prime}\|_\infty^2}{N}
+\frac{\|(1+u^2)G^\prime\|_\infty^2}{N^{1/3}}
\left.\rule{0pt}{12pt}\right)\\
&+&e^T\tilde{\CJ}\left(\rule{0pt}{12pt}\right.
\varepsilon^2 N^{2}\|G^{\prime\prime}\|_\infty^2
+\varepsilon N^{2}\|(1+u^2)G^\prime\|_\infty^2
+\varepsilon^2 N^{4}\|G^\prime\|_1^2
\left.\rule{0pt}{12pt}\right)
\end{eqnarray*}
by Lemma \ref{key result}.
Of course, the same inequality holds if $N$ is replaced by ${N_\delta}$
such that
$$\sum_{i=1}^4
\int_0^\infty\hspace{-5pt}\dd t\,
\ee_\varepsilon\left[
R_{\varepsilon,N_\delta}^{G^\prime,i}(t)
\right]^2
\,\le\,\frac{\delta}{32}\,+\,
e^T\tilde{\CJ}\left(\rule{0pt}{12pt}\right.
\varepsilon^2 N_\delta^{2}\|G^{\prime\prime}\|_\infty^2
+\varepsilon N_\delta^{2}\|(1+u^2)G^\prime\|_\infty^2
+\varepsilon^2 N_\delta^{4}\|G^\prime\|_1^2
\left.\rule{0pt}{12pt}\right)$$
by the choice of $N_\delta$ at the beginning of this proof.
So, choosing $\varepsilon=\varepsilon_{N,{N_\delta}}$
small enough such that both (\ref{epsilon of N and tilde N}) and
$$2\cdot 8e^T\tilde{\CJ}\left(\rule{0pt}{12pt}\right.
\varepsilon^2{N_\delta}^{2}\cdot\|G^{\prime\prime}\|_\infty^2
+\varepsilon{N_\delta}^{2}\cdot c_G\|G^\prime\|_\infty^2
+\varepsilon^2{N_\delta}^{4}\cdot\|G^\prime\|_1^2
\left.\rule{0pt}{12pt}\right)
\,\le\,\delta/4$$
yields
\begin{eqnarray*}
&&\int_0^T\dd t\,\eeo
\left[
\int_{0}^{t}\hspace{-5pt}\int_\NR G^\prime(u)\,
\left(\rule{0pt}{12pt}
(Y_{s}\star J_{{N_\delta}})^2(u)
-
(Y_s\star J_{N})^2(u)\right)
\,\dd u\dd s\right]^2\\
&\le&\rule{0pt}{20pt}
\delta\,+\,
8e^T\tilde{\CJ}\left(\rule{0pt}{12pt}\right.
\frac{\|(1+u^2)G^{\prime\prime\prime}\|_\infty^2}{N^2}
+\frac{\|(1+u^2)G^{\prime\prime}\|_\infty^2}{N}
+\frac{\|(1+u^2)G^\prime\|_\infty^2}{N^{1/3}}
\left.\rule{0pt}{12pt}\right)\\
&\le&\rule{0pt}{20pt}
\delta\,+\,8e^T\tilde{\CJ}\,N^{-1/3}
\sum_{m=1}^3\sup_{u}
|(1+u^2)\frac{\partial^m}{\partial u^m}G(u)|^2.
\end{eqnarray*}
Repeating the above procedure with respect to $N_\delta\ge\tilde{N}$
gives the same inequality for $\tilde{N}$. Hence
\begin{eqnarray*}
&&\int_0^T\dd t\,\eeo
\left[
\int_{0}^{t}\hspace{-5pt}\int_\NR G^\prime(u)\,
\left(\rule{0pt}{12pt}
(Y_{s}\star J_{\tilde{N}})^2(u)
-
(Y_s\star J_{N})^2(u)\right)
\,\dd u\dd s\right]^2\\
&\le&\rule{0pt}{20pt}
4\delta\,+\,32e^T\tilde{\CJ}\,N^{-1/3}
\sum_{m=1}^3\sup_{u}
|(1+u^2)\frac{\partial^m}{\partial u^m}G(u)|^2.
\end{eqnarray*}
for arbitrary but fixed $N,\tilde{N}$ with
$\tilde{N}\ge N$ which finally proves the
lemma since $\delta$ can be made arbitrarily small.\hfill$\Box$

\medskip
\noindent
{\bf Proof of Proposition \ref{cont_mart}(i).}
In this proof the notation $const$ is used when a notation for a
constant is needed thus $const$ can take different values depending on
the situation.

Fix $G\in{\mathscr S}(\NR)$. Applying (\ref{converge}),
there exists a subsequence $(N_k)_{k=1}^\infty$ and 
a measurable subset ${\cal T}_G\subseteq[0,T]$ with $\ell({\cal T}_G)=T$
such that
\begin{equation}\label{conv by t}
\lim_{k\to\infty}\eeo
\left[\rule{0pt}{12pt}\tilde{M}_t^G-\mathfrak{M}_{N_k}(Y)_t^G\right]^2\,=\,0
\end{equation}
for all $t\in{\cal T}_G$. 
For technical reasons assume $T\notin{\cal T}_G$ and let
$\{t_1,t_2,\dots\}\subseteq{\cal T}_G$ be a dense subset of $[0,T]$. 

First observe that $\tilde{M}^G_{t_n}$ 
is ${\cal F}_{t_n}^Y$-\,measurable, $n=1,2,\dots$, and
and the key is to show the following ${\cal F}_t^Y$-\,martingale property
$$\eeo X[\tilde{M}^G_{t_n}-\tilde{M}^G_{t_{n^\prime}}]\,=\,0$$
for $t_{n^\prime},t_n\in\{t_1,t_2,\dots\}$ satisfying
$t_{n^\prime}<t_n$ and an arbitrary random variable $X$ of the form
$X=f(Y_{s_1}(H_1),\dots,Y_{s_p}(H_p))$ 
where $f:\NR^p\to\NR$ is a bounded continuous function, 
$H_i\in{\mathscr S}(\NR)$ and $0\le s_i\le t_{n^\prime},\,1\le i\le p$.
Of course, this martingale property is satisfied if there exists
$const>0$ such that
\begin{equation}\label{mart to show}
\left(\rule{0pt}{12pt}
\eeo X[\tilde{M}^G_{t_n}-\tilde{M}^G_{t_{n^\prime}}]\right)^2
\,\le\,const\cdot\delta
\quad\mbox{for all $\delta>0$.}
\end{equation}

In order to prove (\ref{mart to show}), fix an arbitrary $\delta>0$
and remark that Lemma \ref{key result} implies
$$\int_0^T\dd t\,\ee_\varepsilon\hspace{-3pt}
\left[R_{\varepsilon,N}^{G^\prime,1}(t)\right]^2=O(N^{-1})
\quad\mbox{and}\quad
\int_0^T\dd t\,\ee_\varepsilon\hspace{-3pt}
\left[R_{\varepsilon,N}^{G^\prime,3}(t)\right]^2=O(N^{-1/3})$$
uniformly in $\varepsilon>0$.
Hence, for some $\tau>0$ satisfying $t_n+2\tau<T$,
one can choose $k$ big enough such that both
\begin{equation}\label{trick difference}
\ell(\{t\in[0,T]:
\ee_\varepsilon\hspace{-3pt}
\left[R_{\varepsilon,N_k}^{G^\prime,1}(t)\right]^2
+\;
\ee_\varepsilon\hspace{-3pt}
\left[R_{\varepsilon,N_k}^{G^\prime,3}(t)\right]^2
\ge\delta\})\,\le\,\tau/2
\quad\mbox{for all $\varepsilon>0$}
\end{equation}
and
\begin{equation}\label{delta distance}
\eeo
\left[\tilde{M}^G_{t_n}-\mathfrak{M}_{N_k}(Y)_{t_n}^G\rule{0pt}{12pt}\right]^2
+\;
\eeo
\left[\tilde{M}^G_{t_{n^\prime}}-\mathfrak{M}_{N_k}(Y)_{t_{n^\prime}}^G\right]^2
\,<\,\delta
\end{equation}
hold true. This $k=k_\delta$
is chosen and fixed for proving (\ref{mart to show})
in what follows.

Of course, applying Cauchy-Schwarz, (\ref{delta distance}) implies
\begin{equation}\label{interstep}
\left(\rule{0pt}{12pt}
\eeo X[\tilde{M}^G_{t_n}-\tilde{M}^G_{t_{n^\prime}}]\right)^2
\,\le\,const\left\{\delta+
\left(\rule{0pt}{12pt}
\eeo
X[\mathfrak{M}_{N_k}(Y)_{t_n}^G-\mathfrak{M}_{N_k}(Y)_{t_{n^\prime}}^G]\right)^2
\right\}.
\end{equation}
Now, substituting the definition of $\mathfrak{M}_{N_k}$, one obtains that
\begin{eqnarray*}
\left(\rule{0pt}{12pt}
\eeo
X[\mathfrak{M}_{N_k}(Y)_{t_n}^G-\mathfrak{M}_{N_k}(Y)_{t_{n^\prime}}^G]\right)^2
&=&\left(
\eeo
X[Y_{t_n}(G)-Y_{t_{n^\prime}}(G)]
\;-\int_{t_{n^\prime}}^{t_n}\hspace{-3pt}
\eeo XY_s(G^{\prime\prime})\,\dd s\right.\\
&&\left.
+\;{\gamma}\hspace{-2pt}
\int_{t_{n^\prime}}^{t_n}\hspace{-5pt}\int_\NR G^\prime(u)\,
\eeo X(Y_s\star J_{N_k})^2(u)\,\dd u\dd s
\right)^{\hspace{-2pt}2}
\end{eqnarray*}
where
$$\eeo X(Y_s\star J_{N_k})^2(u)
\,=\,\lim_{\varepsilon\downarrow 0}\,
\hat{\ee}_{\varepsilon}X(Y_s\star J_{N_k})^2(u)$$
such that
$$|\hat{\ee}_{\varepsilon}X(Y_s\star J_{N_k})^2(u)|^2
\,\le\,
\mbox{$\sup_{x\in\NR^p}$}|f(x)|\,
\hat{f}(\|J_{N_k})\|_2^2)$$
for all $\varepsilon\le 1,\,s\in[0,T]$ and $u\in\NR$ 
by Lemma \ref{uni bound} in the Appendix.
Here $f$ is the function defining $X$ while $\hat{f}$ corresponds to
Lemma \ref{uni bound} applied to $(Y_s\star J_{N_k})^2(u)$
and does not depend on $u$. So
$$\int_{t_{n^\prime}}^{t_n}\hspace{-5pt}\int_\NR G^\prime(u)\,
\eeo X(Y_s\star J_{N_k})^2(u)\,\dd u\dd s
\,=\,\lim_{\varepsilon\downarrow 0}
\int_{t_{n^\prime}}^{t_n}\hspace{-5pt}\int_\NR G^\prime(u)\,
\hat{\ee}_{\varepsilon}X(Y_s\star J_{N_k})^2(u)\,\dd u\dd s$$
by dominated convergence and, 
as similar estimates can be obtained for the remaining but easier terms,
one arrives at
$$\left(\rule{0pt}{12pt}
\eeo
X[\mathfrak{M}_{N_k}(Y)_{t_n}^G-\mathfrak{M}_{N_k}(Y)_{t_{n^\prime}}^G]\right)^2$$
$$=\lim_{\varepsilon\downarrow 0}
\left(
\ee_\varepsilon X^\varepsilon\!\left[\rule{0pt}{12pt}\right.
Y_{t_n}^\varepsilon(G)-Y_{t_{n^\prime}}^\varepsilon(G)\,-\hspace{-2pt}
\int_{t_{n^\prime}}^{t_n}\hspace{-3pt}\left\{
Y_s^\varepsilon(G^{\prime\prime})-
{\gamma}\hspace{-2pt}\int_\NR G^\prime(u)\,
(Y_s^\varepsilon\star J_{N_k})^2(u)\,\dd u\right\}\dd s
\left.\rule{0pt}{12pt}\right]
\right)^{\hspace{-2pt}2}$$
$$=\lim_{\varepsilon\downarrow 0}
\left(
\ee_\varepsilon X^\varepsilon\!\left[\rule{0pt}{12pt}\right.
M_{t_n}^{G,\varepsilon}-M_{t_{n^\prime}}^{G,\varepsilon}
+R_\varepsilon^G(t_n)-R_\varepsilon^G(t_{n^\prime})
+{\gamma}\hspace{-2pt}\sum_{i=0}^4
\left(R_{\varepsilon,N_k}^{G^\prime,i}(t_n)-
R_{\varepsilon,N_k}^{G^\prime,i}(t_{n^\prime})\right)
\left.\rule{0pt}{12pt}\right]
\right)^{\hspace{-2pt}2}$$
using (\ref{approxi})\&(\ref{rewrite diff})
for the last equality and writing $X^\varepsilon$
as a substitute for 
$f(Y_{s_1}^\varepsilon(H_1),\dots,Y_{s_p}^\varepsilon(H_p))$.
Notice that 
$\ee_\varepsilon X^\varepsilon
[M_{t_n}^{G,\varepsilon}-M_{t_{n^\prime}}^{G,\varepsilon}]$ 
disappears by the martingale property.
So, if $\varepsilon_0$ is chosen small enough then
\begin{eqnarray}\label{look back to}
&&\left(\rule{0pt}{12pt}
\eeo
X[\mathfrak{M}_{N_k}(Y)_{t_n}^G-\mathfrak{M}_{N_k}(Y)_{t_{n^\prime}}^G]\right)^2
\nonumber\\
&\le&
const\left\{\rule{0pt}{12pt}\right.
\delta\,+\hspace{-7pt}\sum_{t\in\{t_n,t_{n^\prime}\}}\hspace{-5pt}
\left(\ee_{\varepsilon_0}\hspace{-3pt}
\left[R_{\varepsilon_0}^G(t)\right]^2+\,
\sum_{i=1}^4\ee_{\varepsilon_0}\hspace{-3pt}
\left[R_{\varepsilon_0,N_k}^{G^\prime,i}(t)\right]^2\right)
\left.\rule{0pt}{12pt}\hspace{-4pt}\right\}
\end{eqnarray}
by Cauchy-Schwarz. Also, choose $\varepsilon_0$ small enough such that
$$\ee_{\varepsilon_0}\hspace{-3pt}
\left[R_{\varepsilon_0}^G(t_n)\right]^2+\;
\ee_{\varepsilon_0}\hspace{-3pt}
\left[R_{\varepsilon_0}^G(t_{n^\prime})\right]^2
<\,\delta$$
which is possible by (\ref{taylor esti}).
The next lemma provides estimates for the remaining summands. 
\begin{lemma}\label{lemmi}
Fix $0\le i\le4,\,t\in\{t_n,t_{n^\prime}\}$
and $\tau>0$ satisfying $t_n+2\tau<T$. If
$$\ell(\{t\in[0,T]:
\ee_{\varepsilon}\hspace{-3pt}
\left[R_{\varepsilon,N}^{G^\prime,i}(t)\right]^2
\ge\delta\})\,\le\,\tau/2$$
then there exists $\tilde{t}\in[t,t+2\tau]$ such that
$$\ee_\varepsilon\hspace{-3pt}
\left[R_{\varepsilon,N}^{G^\prime,i}(\tilde{t})\right]^2\,<\,\delta
\quad\mbox{and}\quad
\ee_\varepsilon\hspace{-3pt}
\left[R_{\varepsilon,N}^{G^\prime,i}(\tilde{t})-
R_{\varepsilon,N}^{G^\prime,i}(t)\right]^2\,<\,\delta.$$
\end{lemma}
Indeed, observe that if $\tilde{t}\ge t$ then
$$\ee_\varepsilon\hspace{-3pt}
\left[R_{\varepsilon,N}^{G^\prime,i}(\tilde{t})-
R_{\varepsilon,N}^{G^\prime,i}(t)\right]^2
\,=\,
\ee_\varepsilon\hspace{-3pt}
\left[R_{\varepsilon,N}^{G^\prime,i}(\tilde{t}-t)\right]^2$$
by stationarity and the Markov property.
Now assume the contrary of the lemma's assertion, hence
\begin{eqnarray*}
[t,t+2\tau]&\subseteq&
\{\tilde{t}\in[t,t+2\tau]:
\ee_\varepsilon\hspace{-3pt}
\left[R_{\varepsilon,N}^{G^\prime,i}(\tilde{t})\right]^2
\,\ge\,\delta\}\\
&&\cup\;
\{\tilde{t}\in[t,t+2\tau]:
\ee_\varepsilon\hspace{-3pt}
\left[R_{\varepsilon,N}^{G^\prime,i}(\tilde{t})-
R_{\varepsilon,N}^{G^\prime,i}(t)\right]^2
\,\ge\,\delta\}\\
&=&\{\tilde{t}\in[t,t+2\tau]:
\ee_\varepsilon\hspace{-3pt}
\left[R_{\varepsilon,N}^{G^\prime,i}(\tilde{t})\right]^2
\,\ge\,\delta\}\\
&&\cup\;
\{\tilde{t}\in[t,t+2\tau]:
\ee_\varepsilon\hspace{-3pt}
\left[R_{\varepsilon,N}^{G^\prime,i}(\tilde{t}-t)\right]^2
\,\ge\,\delta\}.
\end{eqnarray*}
Thus, as the Lebesgue measures of each of the sets 
on the last equality's right-hand side are bounded by $\tau/2$,
one obtains that $2\tau\le\tau$ which is a contradiction proving the lemma.

Next, for fixed $N_k$, from Lemma \ref{key result} follows that
$$\int_0^T\dd t\,\ee_\varepsilon\hspace{-3pt}
\left[R_{\varepsilon,N_k}^{G^\prime,2}(t)\right]^2=O(\varepsilon)
\quad\mbox{and}\quad
\int_0^T\dd t\,\ee_\varepsilon\hspace{-3pt}
\left[R_{\varepsilon,N_k}^{G^\prime,4}(t)\right]^2=O(\varepsilon^2)$$
and, additionally taking into account (\ref{extra bound}),
one obtains that
$$\ell(\{t\in[0,T]:
\ee_{\varepsilon_1}\hspace{-3pt}
\left[R_{\varepsilon_1,N_k}^{G^\prime,0}(t)\right]^2
+\;
\ee_{\varepsilon_1}\hspace{-3pt}
\left[R_{\varepsilon_1,N_k}^{G^\prime,2}(t)\right]^2
+\;
\ee_{\varepsilon_1}\hspace{-3pt}
\left[R_{\varepsilon_1,N_k}^{G^\prime,4}(t)\right]^2
\ge\delta\})\,\le\,\tau/2$$
for a sufficiently small $\varepsilon_1>0$.
Thus, because (\ref{trick difference}) holds for all $\varepsilon>0$
and so for $\varepsilon_1$ in particular,
one can estimate
\begin{eqnarray*}
\ee_\varepsilon\hspace{-3pt}
\left[R_{\varepsilon,N_k}^{G^\prime,i}({t})\right]^2
&\le&
2\,\ee_\varepsilon\hspace{-3pt}
\left[R_{\varepsilon,N_k}^{G^\prime,i}(\tilde{t})-
R_{\varepsilon,N_k}^{G^\prime,i}(t)\right]^2
+
2\,\ee_\varepsilon\hspace{-3pt}
\left[R_{\varepsilon,N_k}^{G^\prime,i}(\tilde{t})\right]^2\\
&\le&2\delta+2\delta
\end{eqnarray*}
using Lemma \ref{lemmi} for each $i=0,1,2,3,4$ and $t=t_n,t_{n^\prime}$
where $\tilde{t}$ of course depends on 
the chosen $i$ and $t$.
So, when $\varepsilon_0$ in (\ref{look back to}) is replaced by
the minimum of $\varepsilon_0$ and $\varepsilon_1$, 
it follows that
$$\left(\rule{0pt}{12pt}
\eeo
X[\mathfrak{M}_{N_k}(Y)_{t_n}^G-\mathfrak{M}_{N_k}(Y)_{t_{n^\prime}}^G]\right)^2
\,\le\;
const\cdot\delta$$
which, together with (\ref{interstep}),
proves (\ref{mart to show}).
Hence $(\tilde{M}_{s_j}^G)_{j=1}^m$
is an $({\cal F}^Y_{s_j})_{j=1}^m$\,-\,martingale 
for every finite ordered subset $\{s_1,\dots,s_m\}$ 
of $\{t_1,t_2,\dots\}$.

Now, choose arbitrary $s,t\in{\cal T}_G$ and fix $a>0$.
Without restricting the generality one can assume for a moment that
$s,t$ play the role of $t_{n^\prime},t_n$ chosen in the previous part of
this proof.
Combining Chebyshev's inequality and (\ref{delta distance}) yields
\begin{equation}\label{tricky}
\ppo(|\tilde{M}_t^G-\tilde{M}_s^G|>a)
\,\le\,
\frac{const}{a^2}\cdot\delta
+\ppo(|\mathfrak{M}_{N_k}(Y)_t^G-\mathfrak{M}_{N_k}(Y)_s^G|>a/3)
\end{equation}
for the corresponding $k=k_\delta$.
Remark that the set 
$\{|\mathfrak{M}_{N_k}(Y)_t^G-\mathfrak{M}_{N_k}(Y)_s^G|>a/3\}$ 
is open in 
$D([0,T];{\mathscr S}^\prime(\NR))$ with respect to the uniform topology 
and that convergence in $J_1$ 
to elements of $C([0,T];{\mathscr S}^\prime(\NR))$
is equivalent to uniform convergence.
Thus, by Theorem \ref{bigBG}(i), 
the weak convergence of the measures 
$\hat{\pp}_{\hspace{-3pt}\varepsilon},\,\varepsilon\downarrow 0$, implies
$$\ppo(|\mathfrak{M}_{N_k}(Y)_t^G-\mathfrak{M}_{N_k}(Y)_s^G|>a/3)
\,\le\,
\underline{\rm lim}_{\,\varepsilon\downarrow 0}\,
\hat{\pp}_{\hspace{-3pt}\varepsilon}
(|\mathfrak{M}_{N_k}(Y)_t^G-\mathfrak{M}_{N_k}(Y)_s^G|>a/3)$$
where the $lim\,in\!f$ on the right-hand side is equal to
\begin{eqnarray*}
&&\underline{\rm lim}_{\,\varepsilon\downarrow 0}\,
\pp_{\hspace{-3pt}\varepsilon}\hspace{-2pt}
\left(\left|\rule{0pt}{10pt}\right.\!
Y_t^\varepsilon(G)-Y_s^\varepsilon(G)\,-\hspace{-2pt}
\int_{s}^{t}\hspace{-3pt}\left\{
Y_r^\varepsilon(G^{\prime\prime})-
{\gamma}\hspace{-2pt}\int_\NR G^\prime(u)\,
(Y_r^\varepsilon\star J_{N_k})^2(u)\,\dd u\right\}\dd r
\!\left.\rule{0pt}{10pt}\right|>a/3\right)\\
&=&
\underline{\rm lim}_{\,\varepsilon\downarrow 0}\,
\pp_{\hspace{-3pt}\varepsilon}\hspace{-2pt}
\left(\left|\rule{0pt}{10pt}\right.\!
M_{t}^{G,\varepsilon}-M_{s}^{G,\varepsilon}
+R_\varepsilon^G(t)-R_\varepsilon^G(s)
+{\gamma}\hspace{-2pt}\sum_{i=0}^4
\left(R_{\varepsilon,N_k}^{G^\prime,i}(t)-
R_{\varepsilon,N_k}^{G^\prime,i}(s)\right)
\!\left.\rule{0pt}{10pt}\right|>a/3\right)\\
&\le&
\overline{\rm lim}_{\,\varepsilon\downarrow 0}
\left(
\pp_{\hspace{-3pt}\varepsilon}
(|M_{t}^{G,\varepsilon}-M_{s}^{G,\varepsilon}|>a/6)
+
\frac{36}{a^2}
\ee_\varepsilon\hspace{-3pt}
\left[\rule{0pt}{12pt}\right.\hspace{-3pt}
R_\varepsilon^G(t)-R_\varepsilon^G(s)
+{\gamma}\hspace{-2pt}\sum_{i=0}^4
\left(R_{\varepsilon,N_k}^{G^\prime,i}(t)
-R_{\varepsilon,N_k}^{G^\prime,i}(s)\right)
\hspace{-3pt}\left.\rule{0pt}{12pt}\right]^2
\right)
\end{eqnarray*}
where
$$\ee_\varepsilon\hspace{-3pt}
\left[\rule{0pt}{12pt}\right.\hspace{-3pt}
R_\varepsilon^G(t)-R_\varepsilon^G(s)
+{\gamma}\hspace{-2pt}\sum_{i=0}^4
\left(R_{\varepsilon,N_k}^{G^\prime,i}(t)
-R_{\varepsilon,N_k}^{G^\prime,i}(s)\right)
\hspace{-3pt}\left.\rule{0pt}{12pt}\right]^2
\,\le\,const\cdot\delta
\quad\mbox{for all $\varepsilon<\varepsilon_0\wedge\varepsilon_1$}$$
as in the proof of (\ref{mart to show}). 
Using this to estimate the right-hand side of (\ref{tricky})
yields
\begin{equation}\label{from tricky}
\ppo(|\tilde{M}_t^G-\tilde{M}_s^G|>a)
\,\le\, 
\overline{\rm lim}_{\,\varepsilon\downarrow 0}\,
\pp_{\hspace{-3pt}\varepsilon}
(|M_{t}^{G,\varepsilon}-M_{s}^{G,\varepsilon}|>a/6)
\end{equation}
since $\delta$ can be made arbitrarily small.

Now recall that $s,t\in{\cal T}_G$ were arbitrarily chosen and observe that
$$\pp_{\hspace{-3pt}\varepsilon}
(|M_{t}^{G,\varepsilon}-M_{s}^{G,\varepsilon}|>a/6)
\,\le\,
\frac{6^4}{a^4}
\ee_\varepsilon\hspace{-3pt}
\left|\rule{0pt}{10pt}\right.
M_{t}^{G,\varepsilon}-M_{s}^{G,\varepsilon}
\left.\rule{0pt}{10pt}\right|^4
\,\le\,
\frac{6^4C_{4}}{a^4}\,
\ee_\varepsilon\hspace{-3pt}
\left(\rule{0pt}{10pt}\right.
[M^{G,\varepsilon}]_t-[M^{G,\varepsilon}]_s
\left.\rule{0pt}{10pt}\right)^2$$
by first applying Chebyshev's and then Burkholder-Davis-Gundy's inequality
with constant $C_4$.
Furthermore, it is known in this context 
(see \cite{CLO2001} for example) that
$$\ee_\varepsilon\hspace{-3pt}
\left(\rule{0pt}{10pt}\right.
[M^{G,\varepsilon}]_t-[M^{G,\varepsilon}]_s
\left.\rule{0pt}{10pt}\right)^2
\,\le\,
C(T,G)\{\varepsilon^2+(t-s)^2\}.$$
Hence, by (\ref{from tricky}),
there exists $const$ only depending on $T$ and $G$ such that
\begin{equation}\label{kolm centsove cond}
\ppo(|\tilde{M}_t^G-\tilde{M}_s^G|>a)
\,\le\,
const\cdot a^{-4}(t-s)^2
\end{equation}
for all $a>0$ and $s,t\in{\cal T}_G$.

The next step is to construct a continuous process
$(M^G_t)_{t\in[0,T]}$ such that $\tilde{M}_t^G=M^G_t$ 
$\ppo$-a.s.\ for all $t\in{\cal T}_G$. 
But such a construction can be achieved almost the
same way the continuous version of a process is constructed in the
proof of the Kolmogorov-Chentsov theorem (see \cite{KS1991} for
example). As in this proof, it follows from (\ref{kolm centsove cond}) 
that, for a dense subset $D$ of $[0,T]$,
$\{\tilde{M}^G_t(\omega);t\in D\}$ is uniformly continuous in $t$ for
every $\omega\in\Omega^\star$ where $\Omega^\star$ is an event in
${\cal F}^Y_T$ of $\ppo$-measure one.
But in difference to \cite{KS1991}, $D$ should not be  
the set of dyadic rationals in $[0,T]$ but rather an appropriate subset of
the set $\{t_1,t_2,\dots\}$ chosen at the beginning of this proof.
Then one can define $M^G_t(\omega)=0,\,0\le t\le T$, for
$\omega\notin\Omega^\star$ while, for $\omega\in\Omega^\star$,
$M^G_t(\omega)=\tilde{M}^G_t(\omega)$ if $t\in D$ and
$M^G_t(\omega)=\lim_n\tilde{M}^G_{s_n}(\omega)$
for some $(s_n)_{n=1}^\infty\subseteq D$ with $s_n\to t$
if $t\in[0,T]\setminus D$. This gives indeed a continuous process.

To see that $\tilde{M}_t^G=M^G_t$ a.s.\ for all
$t\in{\cal T}_G$ one splits ${\cal T}_G$ into 
$D$ and ${\cal T}_G\setminus D$.
For $t\in D$ one has $\tilde{M}_t^G=M^G_t$ a.s.\ since
$\ppo(\Omega^\star)=1$.
For $t\in{\cal T}_G\setminus D$ 
and $(s_n)_{n=1}^\infty\subseteq D$ with $s_n\to t$
one has
$M^G_t=\lim_n\tilde{M}^G_{s_n}$ a.s.\ by construction
as well as
$\tilde{M}_t^G=\lim_n\tilde{M}^G_{s_n}$ in probability
by (\ref{kolm centsove cond})
which also gives $\tilde{M}_t^G=M^G_t$ a.s.

Realise that, without restricting the generality,
both ${\cal T}_G$ and $D$ and can be chosen to contain
zero as $\mathfrak{M}_N(Y)_0^G=0$ for all $N$ by definition.
Notice that $D\subseteq\{t_1,t_2,\dots\}$ and
$\tilde{M}^G_{t_n}$ is ${\cal F}_{t_n}^Y$-\,measurable for all $n$
and $\Omega^\star\in{\cal F}^Y_T$.
So $M^G_t$ is ${\cal F}_t$\,-\,mesurable for $t\in D$.
Hence $(M^G_t)_{t\in[0,T]}$ is $\NF$-adapted since it is
continuous and $D$ is dense in $[0,T]$.

Finally, the ${\cal F}_t^Y$-\,martingale property of
$\tilde{M}^G_{t_n},\,n=1,2,\dots$, shown by (\ref{mart to show})
implies that $(M_{s_j}^G)_{j=1}^m$
is an $({\cal F}_{s_j})_{j=1}^m$\,-\,martingale 
for every finite ordered subset $\{s_1,\dots,s_m\}$ of $D$.
All these martingales are square integrable because
$\eeo(\tilde{M}^G_{t_n})^2<\infty$ by the choice of $t_n,\,n=1,2,\dots$,
at the beginning of this proof.
Now choose an arbitrary positive $T^\prime<T$.
Then $(M^G_t)_{t\in[0,T^\prime]}$ is a square integrable
$\NF$\,-\,martingale as the limits used to construct this process
can be interchanged with both expectations and conditional expectations
by Doob's maximal inequality for martingales as there must be an
element of $D$ between $T^\prime$ and $T$.\hfill$\Box$\\

\noindent
{\bf Proof of Proposition \ref{cont_mart}(ii).}
Fix $G\in{\mathscr S}(\NR)$. 
Since $(M^G_t)_{t\in[0,T]}$ is a continuous $\NF$-adapted process
it suffices to show that 
for every positive $T^\prime<T$, when restricted to $[0,T^\prime]$,
the process $M^G$ is an  $\NF$-Brownian motion
with variance $2\|G^\prime\|^2_2$. So, in what follows, $T$ is
identified with some positive $T^\prime<T$ to simplify notation.

Obviously, it remains to show that 
$(M^G_t)^2-2\|G^\prime\|^2_2\cdot t,\,t\in[0,T]$, 
is an $\NF$-martingale.
Recalling the construction of $M^G$ 
in the proof of Proposition \ref{cont_mart}(i) above,
the $\NF$-martingale property already follows from
$$\eeo X[(M^G_t)^2-2\|G^\prime\|^2_2\cdot t
-(M^G_{t^\prime})^2+2\|G^\prime\|^2_2\cdot t^\prime\,]\,=\,0$$
for all $t,t^\prime\in D$ such that $t^\prime<t$ and 
$X=f(Y_{s_1}(H_1),\dots,Y_{s_p}(H_p))$
where $f:\NR^p\to\NR$ is a bounded continuous function,
$H_i\in{\mathscr S}(\NR)$ and $0\le s_i\le t^\prime,\,1\le i\le p$.
Again this is verified by showing that 
\begin{equation}\label{to verify}
\left(\rule{0pt}{12pt}
\eeo X[(M^G_t)^2-2\|G^\prime\|^2_2\cdot t
-(M^G_{t^\prime})^2+2\|G^\prime\|^2_2\cdot t^\prime\,]\right)^2
\,\le\,const\cdot\delta
\quad\mbox{for all $\delta>0$}
\end{equation}
for some $const>0$.
So fix $t,t^\prime\in D$ such that $t^\prime<t$ and observe that 
\begin{eqnarray*}
&&\left(\rule{0pt}{12pt}
\eeo X[(M^G_t)^2-2\|G^\prime\|^2_2\cdot t
-(M^G_{t^\prime})^2+2\|G^\prime\|^2_2\cdot t^\prime\,]\right)^2\\
&\le&const\left\{\delta+
\left(\rule{0pt}{12pt}
\eeo
X[(\mathfrak{M}_{N_k}(Y)_{t}^G)^2-(\mathfrak{M}_{N_k}(Y)_{t^\prime}^G)^2
-2\|G^\prime\|^2_2\cdot(t-t^\prime)]\right)^2
\right\}
\end{eqnarray*}
for some $k=k_\delta$ big enough since the inequality
$$\left(\rule{0pt}{12pt}
\eeo[(M^G_t)^2-(\mathfrak{M}_{N_k}(Y)_{t}^G)^2]\right)^2
\,\le\,
2\,\eeo[M^G_t-\mathfrak{M}_{N_k}(Y)_{t}^G]^2
\left(\rule{0pt}{12pt}
\eeo(M^G_t)^2+\eeo(\mathfrak{M}_{N_k}(Y)_{t}^G)^2\right)$$
holds for $t$ and $t^\prime$. Furthermore,
using Lemma \ref{uni bound} in the Appendix
as in the proof of Lemma \ref{cauchy} gives
\begin{eqnarray*}
&&\eeo X(\mathfrak{M}_{N_k}(Y)_{t}^G)^2\\
&=&
\lim_{\varepsilon\downarrow 0}\,
\hat{\ee}_\varepsilon
X\hspace{-3pt}\left(\rule{0pt}{12pt}\right.\hspace{-4pt}
Y_{t}(G)-Y_0(G)
-\int_{0}^{t}\hspace{-3pt}\left\{
Y_s(G^{\prime\prime})-
{\gamma}\hspace{-2pt}\int_\NR G^\prime(u)\,
(Y_s\star J_{N_k})^2(u)\,\dd u\right\}\dd s
\hspace{-4pt}\left.\rule{0pt}{12pt}\right)^2
\end{eqnarray*}
which simplifies to 
$$\hspace{1cm}
=\,\lim_{\varepsilon\downarrow 0}
\ee_\varepsilon
X^\varepsilon\hspace{-3pt}\left(\rule{0pt}{12pt}\right.\hspace{-4pt}
M_{t}^{G,\varepsilon}
+R_\varepsilon^G(t)
+{\gamma}\hspace{-2pt}\sum_{i=0}^4
R_{\varepsilon,N_k}^{G^\prime,i}(t)
\hspace{-4pt}\left.\rule{0pt}{12pt}\right)^2
\;\mbox{with}\;\;
X^\varepsilon\,=\,
f(Y_{s_1}^\varepsilon(H_1),\dots,Y_{s_p}^\varepsilon(H_p)).$$
As the same equality holds for $t^\prime$, one obtains that
\begin{eqnarray*}
&&\left(\rule{0pt}{12pt}
\eeo X[(M^G_t)^2-2\|G^\prime\|^2_2\cdot t
-(M^G_{t^\prime})^2+2\|G^\prime\|^2_2\cdot t^\prime\,]\right)^2\\
&\le&const\left\{\delta+
\left(\rule{0pt}{12pt}
\ee_\varepsilon
X^\varepsilon[(M^{G,\varepsilon}_t)^2
-(M^{G,\varepsilon}_{t^\prime})^2
-2\|G^\prime\|^2_2\cdot(t-t^\prime)]\right)^2
\right\}
\end{eqnarray*}
for a sufficiently small $\varepsilon>0$ by estimating
$$\ee_\varepsilon
M^{G,\varepsilon}_t
\hspace{-3pt}\left(\rule{0pt}{12pt}\right.\hspace{-4pt}
R_\varepsilon^G(t)
+{\gamma}\hspace{-2pt}\sum_{i=0}^4
R_{\varepsilon,N_k}^{G^\prime,i}(t)
\hspace{-4pt}\left.\rule{0pt}{12pt}\right)
\quad\mbox{and}\quad
\ee_\varepsilon
\hspace{-3pt}\left(\rule{0pt}{12pt}\right.\hspace{-4pt}
R_\varepsilon^G(t)
+{\gamma}\hspace{-2pt}\sum_{i=0}^4
R_{\varepsilon,N_k}^{G^\prime,i}(t)
\hspace{-4pt}\left.\rule{0pt}{12pt}\right)^2$$
for $t$ and $t^\prime$ 
using the bounds derived
in the proof of Proposition \ref{cont_mart}(i).

Now $(M^{G,\varepsilon}_t)^2,\,t\ge 0$,
is a submartingale in the class $(\mathfrak{D}L)$.
Hence
$(M^{G,\varepsilon}_t)^2-\langle M^{G,\varepsilon}\rangle_t,\,t\ge 0$,
is a martingale so that
\begin{eqnarray*}
&&\left(\rule{0pt}{12pt}
\eeo X[(M^G_t)^2-2\|G^\prime\|^2_2\cdot t
-(M^G_{t^\prime})^2+2\|G^\prime\|^2_2\cdot t^\prime\,]\right)^2\\
&\le&const\left\{\delta+
\left(\rule{0pt}{12pt}
\ee_\varepsilon
X^\varepsilon[\langle M^{G,\varepsilon}\rangle_t
-\langle M^{G,\varepsilon}\rangle_{t^\prime}
-2\|G^\prime\|^2_2\cdot(t-t^\prime)]\right)^2
\right\}.
\end{eqnarray*}
Finally 
$\ee_\varepsilon
[\langle M^{G,\varepsilon}\rangle_t
-\langle M^{G,\varepsilon}\rangle_{t^\prime}
-2\|G^\prime\|^2_2\cdot(t-t^\prime)]^2$
can be made arbitrarily small by choosing a suitable $\varepsilon$
which proves (\ref{to verify}) 
hence part (ii) of Proposition \ref{cont_mart}.
The last argument is standard 
and can be found in \cite{CLO2001}, for example.\hfill$\Box$\\

\noindent
{\bf Proof of Proposition \ref{cont_mart}(iii).}
Fix $a_1,a_2\in\NR$ and $G_1,G_2\in{\mathscr S}(\NR)$.
The wanted linearity holds for ${\mathfrak M}_N(Y)$ and,
because ${\mathfrak M}_N(Y)$ is an approximation for
$(\tilde{M}^G)_{G\in{\mathscr S}(\NR)}$,
the linearity should also hold for the version 
$({M}^G)_{G\in{\mathscr S}(\NR)}$ of $(\tilde{M}^G)_{G\in{\mathscr S}(\NR)}$.
But some care has to be taken since the construction of 
$({M}^G)_{G\in{\mathscr S}(\NR)}$ depends on the choice of subsequences
and, also, since the notion of version used in this paper is special 
as not all $t\in[0,T]$ are covered. 

By Proposition \ref{cont_mart}(i),
there are sets 
${\cal T}_{G_1},{\cal T}_{G_2},{\cal T}_{a_1G_1+a_2G_2}$
corresponding to the processes
$M^{G_1}$, $M^{G_2}$, $M^{a_1 G_1+a_2 G_2}$.
First one wants to find a set
$${\cal T}\,\subseteq\,
{\cal T}_{G_1}\cap{\cal T}_{G_2}\cap{\cal T}_{a_1G_1+a_2G_2}
\quad\mbox{dense in $[0,T]$}$$
such that
\begin{equation}\label{firstly}
\tilde{M}_t^{a_1 G_1+a_2 G_2}\,=\,a_1\tilde{M}_t^{G_1}+a_2\tilde{M}_t^{G_2}
\quad\mbox{a.s.}\quad\mbox{for $t\in{\cal T}$.}
\end{equation}
This is achieved by successively choosing subsequences as follows.
Using (\ref{conv by t}), there is a subsequence 
$(k_j)_{j=1}^\infty$ of $(N_k)_{k=1}^\infty$ such that
\begin{equation}\label{secoli}
\tilde{M}_t^{a_1 G_1+a_2 G_2}\,=\,
\lim_{j\to\infty}\left(\rule{0pt}{12pt}\right.
a_1{\mathfrak M}_{k_j}(Y)_t^{G_1}
+a_2{\mathfrak M}_{k_j}(Y)_t^{G_2}
\left.\rule{0pt}{12pt}\right)
\quad\mbox{a.s.}\quad\mbox{for $t\in{\cal T}_{a_1G_1+a_2G_2}$.}
\end{equation}
Now, using (\ref{converge}) with respect to 
$(k_j)_{j=1}^\infty$ and $G_1$,
there is a measurable subset ${\cal T}^\prime_{G_1}\subseteq[0,T]$ with
$\ell({\cal T}^\prime_{G_1})=T$ and a subsequence
$(j_l)_{l=1}^\infty$ of $(k_j)_{j=1}^\infty$ such that
$$\tilde{M}_t^{G_1}\,=\,
\lim_{l\to\infty}{\mathfrak M}_{j_l}(Y)_t^{G_1}
\quad\mbox{a.s.}\quad\mbox{for $t\in{\cal T}^\prime_{G_1}$.}$$
Notice that ${\cal T}^\prime_{G_1}$ and ${\cal T}_{G_1}$
can be different. Similarly, one obtains that
$$\tilde{M}_t^{G_2}\,=\,
\lim_{m\to\infty}{\mathfrak M}_{l_m}(Y)_t^{G_2}
\quad\mbox{a.s.}\quad\mbox{for $t\in{\cal T}^\prime_{G_2}$}$$
where $(l_m)_{m=1}^\infty$ is a subsequence of
$(j_l)_{l=1}^\infty$ and $\ell({\cal T}^\prime_{G_2})=T$.
Then
$${\cal T}\,\stackrel{\mbox{\tiny def}}{=}\,
{\cal T}_{G_1}\cap{\cal T}_{G_2}\cap{\cal T}_{a_1G_1+a_2G_2}
\cap{\cal T}^\prime_{G_1}\cap{\cal T}^\prime_{G_2}
\,\subseteq\,
{\cal T}_{G_1}\cap{\cal T}_{G_2}\cap{\cal T}_{a_1G_1+a_2G_2}$$
and ${\cal T}$ is dense in $[0,T]$ because $\ell({\cal T})=T$.
Furthermore, using the subsequence
$(l_m)_{m=1}^\infty$ instead of $(k_j)_{j=1}^\infty$
in (\ref{secoli}) implies (\ref{firstly}).

But, by Proposition \ref{cont_mart}(i),
(\ref{firstly}) is equivalent to
$$M_t^{a_1 G_1+a_2 G_2}\,=\,a_1 M_t^{G_1}+a_2 M_t^{G_2}
\quad\mbox{a.s.}\quad\mbox{for $t\in{\cal T}$}$$
which proves part (iii) of Proposition \ref{cont_mart}
because the processes $M^{a_1 G_1+a_2 G_2},\,M^{G_1},\,M^{G_2}$
are continuous.\hfill$\Box$\\

\noindent
{\bf Proof of Proposition \ref{cont_mart}(iv).}
Remark that part (iv) would not follow from part (ii) allone but,
including part (iii), it is straight forward to check both the
Gaussian distribution and the covariance structure of the
process  $M_t^G$ indexed by $t\in[0,T]$ and $G\in{\mathscr S}(\NR)$. 
Of course, from the covariance structure follows that the index set of
the process can be extended to $t\in[0,T]$ and absolutely continuous
functions $G$ on $\NR$ with density $G^\prime\in L^2(\NR)$ without
changing the underlying probability space. Hence
$$\tilde{B}(t,u)\,=\,
M_t^{G_u}/\sqrt{2}\,,\quad t\in[0,T],\;u\in\NR,$$
is properly defined using test functions 
$G_u(\tilde{u}),\,\tilde{u}\in\NR$, given by
$$G_u(\tilde{u})\,=\left\{\begin{array}{rcl}
0\vee(u\wedge\tilde{u})&:&u\ge 0,\\
0\wedge(u\vee\tilde{u})&:&u<0.
\end{array}\right.$$
Obviously, $\tilde{B}(t,u),\,t\in[0,T],\,u\in\NR$,
is a centred Gaussian process on 
$(D([0,T];{\mathscr S}^\prime(\NR)),{\cal F}^Y_T,$ $\ppo)$
with covariance 
$\eeo\tilde{B}(t,u)\tilde{B}(t^\prime,u^\prime)
=(t\wedge t^\prime)(|u|\wedge|u^\prime|)$ if $u,u^\prime$ have the
same sign and vanishing covariance otherwise.
So, as in the proof of the  Kolmogorov-Chentsov theorem, one can
construct a version $B(t,u)$ of $\tilde{B}(t,u)$ on the same
probability space which is continuous in $t$ and $u$, hence, 
is a Brownian sheet. 
By standard theory on random linear functionals, 
see \cite{W1986} for a good reference, 
there is an ${\mathscr S}^\prime(\NR)$-valued version of the process 
$M_t^G$ which is of course indistinguishable of
$$\sqrt{2}\int_\NR B(t,u)G^{\prime\prime}(u)\,\dd u
\quad t\in[0,T],\;G\in{\mathscr S}(\NR),$$
finally proving part (iv) of Proposition \ref{cont_mart}.\hfill$\Box$
%
\section{Appendix}
Recall that $\hat{\pp}_{\hspace{-3pt}\varepsilon}$ is the push forward
of $\ppe$ with respect to the map $Y^\varepsilon$ 
introduced on page \pageref{bigBG} and denote by
$\hat{\ee}_{\varepsilon}$ the expectation when integrating against
$\hat{\pp}_{\hspace{-3pt}\varepsilon}$.
Then it is a consequence of Theorem \ref{bigBG}(i)
that weak convergence implies
\begin{equation}\label{still conv weak}
\hat{\ee}_{\varepsilon}X
\;\to\;\eeo X,
\quad\varepsilon\downarrow 0\,,
\end{equation}
for $X=f(Y_{s_1}(H_1),\dots,Y_{s_p}(H_p))$
defined by bounded continuous maps $f:\NR^p\to\NR$ and 
$H_i\in{\mathscr S}(\NR),\,0\le s_i\le T,\,1\le i\le p$, 
although such functions $X$ are not $J_1$\,-\,continuous 
on the space $D([0,T];{\mathscr S}^\prime(\NR))$.

The lemma below states that the boundedness condition on $f$ can be
relaxed when the one-dimensional marginals of the limit process are
Gaussian. This result is not new but the specific statement needed in
this paper could not be found in the literature. Remark that if the
limit process does not have Gaussian one-dimensional marginals then,
for polynomial singularities,
instead of weak convergence of measures one should consider
convergence in Wasserstein spaces.
\begin{lemma}\label{uni bound}
The convergence (\ref{still conv weak}) remains true 
for $X$ defined by continuous functions $f$ with polynomial growth
and 
$${\rm sup}_{\varepsilon\le 1}
|\hat{\ee}_{\varepsilon}X|^2
+
|\eeo X|^2
\,\le\,\hat{f}(\|H_1\|_2^2\,,\dots,\|H_p\|_2^2)$$
where $\hat{f}$ is a polynomial not depending on the time points
$s_1,\dots,s_p$ defining $X$.
\end{lemma}
{\it Proof}. 
It suffices to show the lemma for polynomials $f$. The convergence
claim follows from Theorem \ref{bigBG}(ii). Indeed, as the
one-dimensional marginal distributions of $Y$ under $\ppo$ are
Gaussian, one can cut-off $f$ turning it into a bounded continuous
function for which (\ref{still conv weak}) holds and estimate the remainder using the
exponential decay of the tails of the Gaussian distribution.

The uniform bound $\hat{f}(\|H_1\|_2^2\,,\dots,\|H_p\|_2^2)$ also
follows from Theorem \ref{bigBG}(ii) by successively applying
H\"older's inequality and estimating moments of Gaussian distributions
by powers of the variances.
Notice that the supremum is taken over $0<\varepsilon\le 1$
but any other bounded subset of $\varepsilon>0$ 
could have been used.\hfill$\Box$

\newpage

\end{document}